
\documentstyle{amsppt}
\NoBlackBoxes \magnification1200

\pageheight{9 true in} \pagewidth{6.5 true in}
\def\Lap{{\Cal L}}
\def\vtheta{\vartheta}
\def\vphi{\varphi}

\define\sumstar{\sideset \and^* \to\sum}
\define\sumflat{\sideset \and^{\flat} \to\sum}

\topmatter
\title
An uncertainty principle for arithmetic sequences
\endtitle
\author
Andrew Granville and K. Soundararajan
\endauthor
\address{D{\'e}partment  de
Math{\'e}matiques et Statistique, Universit{\'e}
de Montr{\'e}al, CP 6128 succ
Centre-Ville, Montr{\'e}al, QC  H3C 3J7, Canada}\endaddress
\email{andrew{\@}dms.umontreal.ca}
\endemail
\address{Department of Mathematics, University of Michigan, Ann Arbor,
Michigan 48109, USA} \endaddress \email{ksound{\@}umich.edu}
\endemail

\abstract Analytic number theorists usually seek to show that
sequences which appear naturally in arithmetic are
``well-distributed'' in some appropriate sense.  In various
discrepancy problems, combinatorics researchers have analyzed
limitations to equi-distribution, as have Fourier analysts when
working with the ``uncertainty principle''. In this article we
find that these ideas have a natural setting in the analysis of
distributions of sequences in analytic number theory, formulating
a general principle, and giving several examples.
\endabstract

\thanks{Le premier auteur est partiellement soutenu par une bourse
du Conseil  de recherches en sciences naturelles et en g\' enie du
Canada. The second  author is partially supported by the National
Science Foundation.}
\endthanks
\endtopmatter

\document

\head 1. Introduction \endhead

\noindent In this paper we investigate the limitations to the
equidistribution of interesting ``arithmetic sequences'' in
arithmetic progressions and short intervals.  Our discussions are
motivated by a general result of K.F. Roth [15] on irregularities of
distribution, and a particular result of H. Maier [11] which imposes
restrictions on the equidistribution of primes.

If ${\Cal A}$ is a subset of the integers in $[1,x]$ with $|{\Cal
A}| =\rho x$  then Roth proved that there exists $N\le x$ and an
arithmetic progression $a \pmod q$ with $q\le \sqrt{x}$ such that
$$
\Big| \sum\Sb n\in {\Cal A},\ n\leq N\\ n\equiv a\pmod q\endSb 1
-\frac{1}{q} \sum\Sb n\in {\Cal A}\\ n\leq N\endSb 1\ \Big| \gg
\sqrt{\rho(1-\rho)} x^{\frac 14}.
$$
In other words, keeping away from sets of density $0$ or $1$,
there must be an arithmetic progression in which the number of
elements of ${\Cal A}$ is a little different from the average.
Following work of A. Sarkozy and J. Beck,  J. Matousek and J. Spencer
[12] showed that Roth's theorem is best possible, in that there is
a set $\Cal A$ containing $\sim x/2$ integers up to $x$, for which
$| \# \{n\in {\Cal A}:\ n\leq N,\  n\equiv a\pmod q\}  - \# \{n\in
{\Cal A}:\ n\leq N\}/q | \ll x^{1/4}$ for all $q$ and $a$ with
$N\leq x$.

Roth's result concerns arbitrary sequences of integers, as
considered in combinatorial number theory and harmonic analysis.
We are more interested here in sets of integers that arise in
arithmetic, such as the primes.  In [11] H. Maier developed an
ingenious method to show that for any $A\ge 1$ there are
arbitrarily large $x$ such that the interval $(x,x+(\log x)^A)$
contains significantly more primes than usual (that is, $\ge
(1+\delta_A) (\log x)^{A-1}$ primes for some $\delta_A >0$) and
also intervals $(x,x+(\log x)^{A})$ containing significantly fewer
primes than usual.   Adapting his method J. Friedlander and A.
Granville [3] showed that there are arithmetic progressions containing
significantly more (and others with significantly fewer) primes
than usual.  A weak form of their result is that, for every $A\ge
1$ there exist large $x$ and an arithmetic progression $a\pmod q$
with $(a,q)=1$ and $q\le x/(\log x)^A$ such that
$$
\Big| \pi(x;q,a) -\frac{\pi(x)}{\phi(q)} \Big| \gg_A
\frac{\pi(x)}{\phi(q)}. \tag{1.1}
$$
If we compare this to Roth's bound   we note two differences: the
discrepancy exhibited is much larger in (1.1) (being within a
constant factor of the main term), but the modulus $q$ is much
closer to $x$ (but not so close as to be trivial).

Recently A. Balog and T. Wooley [1] proved that the sequence of
integers that may be written as the sum of two squares also
exhibits ``Maier type'' irregularities in some intervals
$(x,x+(\log x)^A)$ for any fixed, positive $A$. While previously
Maier's results on primes had seemed inextricably linked to the
mysteries of the primes, Balog and Wooley's example suggests that
such results should be part of a general phenomenon.   Indeed, we
will provide here a general framework for such results on
irregularities of distribution, which will include, among other
examples, the sequence of primes and the sequence of sums of two
squares.  Our results may be viewed as an ``uncertainty
principle'' which establishes that most arithmetic sequences of
interest are either not-so-well distributed in longish arithmetic
progressions, or are not-so-well distributed in both short
intervals and short  arithmetic progressions.

\subhead 1a. Examples \endsubhead  

We now highlight this
phenomenom with several examples:\ For a given set of integers
$\Cal A$, let $\Cal A(N)$ denote the number of elements of $\Cal
A$ which are $\leq N$, and $\Cal A(N;q,a)$ denote those that are
$\leq N$ and $\equiv a \pmod q$.

\bigskip

$\bullet$\  We saw in Maier's theorem that the primes are not so
well-distributed. We might ask whether there are subsets
$\Cal A$ of the primes up to $x$ which are well-distributed? Fix
$u\geq 1$. We show that for any $x$ there exists $y\in (x/4,x)$
such that either
$$|\Cal A(y)/y-\Cal A(x)/x|\gg_u \Cal A(x)/x \tag{1.2a} $$
(meaning that the subset is poorly distributed in short
intervals), or there exists some arithmetic progression $a \pmod
\ell$ with $(a,\ell)=1$ and $\ell \le x/(\log x)^u$, for which
$$
\Big| {\Cal A}(y;\ell,a)- \frac{\Cal A(y)}{\phi(\ell)}\Big| \gg_u
\frac{{\Cal A}(x)}{\phi(\ell)}. \tag{1.2b}
$$
In other words, we find ``Maier type'' irregularities in the
distribution of {\sl any subset of the primes}. (If we had chosen
${\Cal A}$ to be the primes $\equiv 5 \pmod 7$ then this is of no
interest when we take $a=1, \ell=7$. To avoid this minor
technicality we can add ``For a given finite set of ``bad primes''
${\Cal S}$, we can choose such an $\ell$ for which $(\ell,\Cal
S)=1$''. Here and henceforth $(\ell,\Cal S)=1$ means that
$(\ell,p)=1$ for all $p\in \Cal S$. )

\bigskip

$\bullet$\ With probability 1 there are {\sl no} ``Maier type''
irregularities in the distribution of randomly chosen subsets of
the integers. Indeed such irregularities seem to depend on the
subset having some arithmetic structure.  So instead of taking
subsets of all the integers, we need to take subsets of a set
which already has some arithmetic structure. For example, define
$\Cal S_\epsilon$ to be the set of integers $n$ having no prime
factors in the interval $[(\log n)^{1-\epsilon}, \log n]$, so that
$\Cal S_\epsilon(N)\sim (1-\epsilon)N$. Notice that the primes are
a subset of  $\Cal S_\epsilon$.  Our results imply that any subset
$\Cal A$ of $\Cal S_\epsilon$ is poorly distributed in that for
any $x$ there exists $y\in (x/4,x)$ such that either (1.2a) holds,
or there exists some arithmetic progression $a \pmod \ell$   and
$\ell \le x/(\log x)^u$ with $(a,\ell)=1$, for which a suitably
modified (1.2b) holds (that is with $\phi(\ell)$ replaced by $\ell
\prod_{p|\ell,\ (\log x)^{1-\epsilon}<p< \log x} (1-1/p)$).

\bigskip

$\bullet$\ Let $K$ be an algebraic number field with $[K:{\Bbb Q}]>1$.  
Let $R$ denote the ring of integers of $K$ and let $C$ be an 
ideal class from the class group of $R$. Take $\Cal A$ be the
set of positive integers which are the norm of some (integral) ideal 
belonging to $C$.  (In Balog and Wooley's example, 
$\Cal A$ is the set of numbers of the form $x^2+y^2$, 
with $C$ the class of principal ideals in $R=\Bbb Z[i]$.) 
From our work it follows that the set $\Cal A$
is poorly distributed in arithmetic progressions; that is, a
suitably modified version of (1.2b) holds.
Moreover, 
if we replace $R$ by any order in $K$ then either (1.2a) holds 
or a suitably modified version of (1.2b) holds (and we expect 
that, with some effort, one can prove that the 
suitably modified (1.2b) holds).

\bigskip

$\bullet$\ Let $\Cal B$ be a given set of $x$ integers and $\Cal
P$ be a given set of primes. Define $\Cal S(\Cal B,\Cal P,z)$ to
be the number of integers in $\Cal B$ which do not have a prime
factor $p\in \Cal P$ with $p\leq z$. Sieve theory is concerned
with estimating $\Cal S(\Cal B,\Cal P,z)$ under certain natural
hypothesis  for $\Cal B, \Cal P$ and $u:=\log x/\log z$. The
fundamental lemma of sieve theory (see [7]) implies (for example
when $\Cal B$ is the set of integers in an interval) that
$$
\left| \Cal S(\Cal B,\Cal P,z) - x \prod\Sb p\in \Cal P, p\leq z\endSb
\left( 1 -\frac 1p \right)\right|
\ll \left( \frac{1+o(1)}{u\log u} \right)^u x \prod\Sb p\in \Cal P, p\leq z\endSb
\left( 1 -\frac 1p \right)
$$
for $u<z^{1/2+o(1)}$.
It is known that this result is essentially ``best-possible'' in that one can
construct examples for which the bound is obtained (both as an 
upper and lower bound).
However these bounds are obtained in quite special examples, and one might
suspect that in many cases which one encounters, those bounds might be
significantly sharpened. It turns out that these bounds cannot be 
improved for intervals
$\Cal B$, when $\Cal P$ contains at least a positive proportion of the primes:

\proclaim{Corollary 1.1} Suppose that $\Cal P$ is a given set of
primes for which $\# \{p\in \Cal P:\ p\leq y\} \gg \pi(y)$ for all
$y\in (\sqrt{z},z]$. There exist constants $c>0$ such that for any
$u\ll \sqrt{z}$ there exist intervals $I_\pm$ of length $\geq z^u$
for which
$$
\align \Cal S(I_+,\Cal P,z) &\geq \left\{ 1 + \left(
\frac{c}{u\log u} \right)^u\right\}
|I_+| \prod\Sb p\in \Cal P, p\leq z\endSb   \left( 1 -\frac 1p \right) \\
\text{\rm and} \quad \Cal S(I_-,\Cal P,z) &\leq \left\{ 1 -
\left( \frac{c}{u\log u} \right)^u\right\} |I_-| \prod\Sb p\in
\Cal P, p\leq z\endSb   \left( 1 -\frac 1p \right) .
\endalign
$$
Moreover if $u\leq (1-o(1))\log\log z/\log\log\log z$ then our
intervals $I_\pm$ have length $\leq z^{u+2}$.
\endproclaim

\bigskip

$\bullet$\ What about sieve questions in which the set of
primes does not have positive lower density (in the set of
primes)? If $\Cal P$ contains too few primes then we should expect
the sieve estimate to be very accurate, so we must insist on some
lower bound: for instance that if $q=\prod_{p\in \Cal P} p$ then
$$
\sum_{p|q} \frac{\log p}{p} \ge 60 \log \log \log q. \tag{1.3}
$$
(Note that $\sum_{p|q} (\log p)/p \le (1+o(1))\log \log q$, the
bound being attained when $q$ is the product of the primes up to
some large $y$.)

\proclaim{Corollary 1.2} Let $q$ be a large square-free number,
which satisfies {\rm (1.3)}, and define $z:= (\prod_{p|q}
p^{1/p})^{c_1}$ for a certain constant $c_1>0$ . There exists a
constant $c_2>0$ such that if $\sqrt{z} \ge u\gg (\log\log q/\log
z)^3$ then exist intervals $I_\pm$ of length at least $z^u$ such
that
$$
  \sum\Sb n\in I_+\\ (n,q)=1 \endSb 1 \ge   \{ 1+ 1/u^{c_2u}
\} \frac{\phi(q)}{q} |I_+| ,\quad \text{\rm and} \quad  \sum\Sb
n\in I_-\\ (n,q)=1 \endSb 1 \le
 \{ 1- 1/u^{c_2u} \} \frac{\phi(q)}{q} |I_-| .
$$
\endproclaim

\bigskip

$\bullet$\ The  reduced residues $\pmod q$ are expected to be
distributed much like random numbers chosen with probability
$\phi(q)/q$. Indeed when $\phi(q)/q\to 0$ this follows from work
of C. Hooley [10]; and  of H.L. Montgomery and R.C. Vaughan [13]
who showed that $\# \{ n\in [m,m+h):\ (n,q)=1\}$ has  Gaussian
distribution with mean and variance equal to $h\phi(q)/q$, as $m$
varies over the integers, provided $h$ is suitably large. This
suggests that $\# \{ n\in [m,m+h):\ (n,q)=1\}$ should be
$\{1+o(1)\}(h\phi(q)/q)$ provided $h\ge \log^2q$, however by
Corollary 1.2 this is not true for $h=\log^Aq$ for any given $A>0$,
provided that $\sum_{p|q} (\log p)/p \gg \log\log q$ (a condition
satisfied by many highly composite $q$).
\bigskip

In \S 6 we shall give further new examples of sequences to which
our results apply.

\subhead 1b. General results \endsubhead

Our main result (Theorem 3.1) is too technical to introduce at
this stage. Instead we motivate our setup (postponing complete
details to \S 2) and explain some consequences.

Let ${\Cal A}$ denote a sequence $a(n)$ of non-negative real
numbers. We are interested in determining whether the $a(n)$ are
well-distributed in short intervals and in arithmetic
progressions, so let  ${\Cal A}(x) =\sum_{n\le x} a(n)$ (so if
${\Cal A}$ is a set of positive integers  then $a(n)$ is its
indicator function). Thinking of ${\Cal A}(x)/x$ as the average
value of $a(n)$, we may expect that if  ${\Cal A}$ is
well-distributed in short intervals then
$$
{\Cal A}(x+y) - {\Cal A}(x) \approx y\frac{{\Cal A}(x)}{x},
\tag{1.4}
$$
for suitable $y$.

To understand the distribution of ${\Cal A}$ in arithmetic
progressions, we begin with those $n$ divisible by $d$. We will
 suppose that the proportion of ${\Cal A}$ which is
divisible by $d$ is approximately $h(d)/d$ where $h(.)$ is a
non-negative multiplicative function; in other words,
$$
{\Cal A}_d(x) := \sum\Sb n\le x\\ d|n \endSb a(n) \approx
\frac{h(d)}{d} \Cal A(x),  \tag{1.5}
$$
for each $d$ (or perhaps when $(d,{\Cal S})=1$, where ${\Cal S}$
is a finite set of `bad' primes). The reason for taking $h(d)$ to
be a multiplicative function is that for most sequences that
appear in arithmetic one expects that the criterion 
of being divisible by an integer $d_1$ should be 
``independent'' of the criterion of being divisible 
by an integer $d_2$ coprime to $d_1$.

If the asymptotic behavior
 of ${\Cal A}(x;q,a)$ for $(q,{\Cal S})=1$
depends only on the g.c.d. of $a$ and $q$ then, by (1.5), we
arrive at the prediction that, for $(q,{\Cal S})=1$,
$$
{\Cal A}(x;q,a) \approx \frac{f_q(a)}{q\gamma_q} {\Cal A}(x), \tag{1.6}
$$
where $\gamma_q = \prod_{p|q} ((p-1)/(p-h(p)))$ and
$f_q(a)$ is a certain non-negative multiplicative function of $a$
for which $f_q(a)=f_q((a,q))$ (thus $f_q(a)$ is periodic $\pmod q$).  
In \S 2 we shall give an explicit
description of $f_q$ in terms of $h$.

In the spirit of Roth's theorem we ask how good is the
approximation (1.6)?  And, in the spirit of Maier's theorem we ask
how good is the approximation  (1.4)?
\medskip

\noindent{\bf Example 1.} We take $a(n)=1$ for all $n$.  We may take ${\Cal S}
=\emptyset$ and $h(n)=1$ for all $n$.  Then $f_q(a)=1$ for all $q$ and all $a$,
and $\gamma_q=1$.  Clearly both (1.6) and (1.4) are good approximations
with an error of at most $1$.

\noindent{\bf Example 2.}  We take $a(n)=1$ if $n$ is prime and
$a(n)=0$ otherwise. Then we may take ${\Cal S}=\emptyset$ and
$h(n)=1$ if $n=1$ and $h(n)=0$ if $n>1$.  Further $f_q(a)=1$ if
$(a,q)=1$ and $f_q(a)=0$ otherwise, and $\gamma_q = \phi(q)/q$.
The approximation (1.6) is then the prime number theorem for
arithmetic progressions for small $q\le (\log x)^A$.  Friedlander
and Granville's result (1.1) sets limitations to (1.6), and
Maier's result sets limitations to (1.4).

\noindent {\bf Example 3.} Take $a(n)=1$ if $n$ is the sum of
two squares and $a(n)=0$ otherwise.  Here we take ${\Cal S}=\{2\}$,
and for odd prime powers $p^k$ we have $h(p^k)=1$ if $p^k\equiv
1\pmod 4$ and $h(p^k)=1/p$ otherwise.  Balog and Wooley's result
places restrictions on the validity of (1.4).


\proclaim{Corollary 1.3}  Let
$\Cal A$, ${\Cal S}$, $h$, $f_q$ and $\gamma_q$ be as above.
Let $x$ be sufficiently large and in particular suppose that
${\Cal S}\subset [1,\log \log x]$.
Suppose that $0\le h(n)\le 1$ for all $n$.  Suppose that
$$
\sum_{p\leq \log x} \frac{1-h(p)}{p} \log p \ge \alpha \log \log x, \tag{1.7}
$$
for some $\alpha \ge 60\log \log \log x/\log \log x$ and set $\eta =\min(\alpha/3,1/100)$.
Then for each $5/\eta^2 \le u\le \eta (\log x)^{\eta/2}$
there exists $y\in (x/4,x)$ and an arithmetic progression
$a\pmod{\ell}$ with $\ell \le {x}/{(\log x)^u}$ and $(\ell,{\Cal S})=1$
such that
$$
\Big| {\Cal A}(y;\ell, a) - \frac{f_\ell(a)}{\ell\gamma_\ell}
y\frac{{\Cal A}(x)}{x}\Big|
\gg \exp\Big( -\frac{u}{\eta} (1+25\eta) \log (2u/\eta^3) \Big)
\frac{{\Cal A}(x)}{\phi(\ell)}.
$$
\endproclaim

\remark{Remarks} Since the Corollary appears quite technical, some
explanation is in order.

\smallskip
$\bullet$ The condition $0\le h(n)\le 1$ is not as restrictive as it might 
appear.
We will show in Proposition 2.1 if there are many
primes with $h(p)>1$ then it is quite easy to construct
large discrepancies for the sequence ${\Cal A}$.
\smallskip
$\bullet$ The condition (1.7) ensures that $h(p)$ is not always
close to $1$; this is essential in order to eliminate
the very well behaved Example 1.
\smallskip
$\bullet$  The conclusion of the
Corollary may be weakly (but perhaps more transparently)
written as
$$
\Big|{\Cal A}(y;\ell,a) -\frac{f_\ell(a)}{\ell \gamma_\ell}
y\frac{{\Cal A}(x)}{x}
\Big|  \gg_{\alpha, u} \frac{{\Cal A}(x)}{\phi(\ell)}.
$$
\smallskip
$\bullet$  The lower bound given is a multiple of ${\Cal
A}(x)/\phi(\ell)$, rather than of the main term
$(f_\ell(a)/\ell\gamma_\ell)(y{\Cal A}(x)/x)$. The main reason for
this is that $f_\ell(a)$ may well be $0$, in which case such a
bound would have no content. In fact, since $(y/x)<1$ and
$\phi(\ell)\leq \ell\gamma_\ell$, so the function used is larger
and more meaningful than  the main term itself.
\smallskip
$\bullet$  It might appear more natural to compare
${\Cal A}(y;\ell,a)$ with $(f_\ell(a)/\ell\gamma_\ell){\Cal A}(y)$.
In most examples that we consider the average ${\Cal A}(x)/x$ ``varies
slowly'' with $x$, so we expect little difference between ${\Cal A}(y)$
and $y{\Cal A}(x)/x$ (we have   $\sim 1/\log x$ in Example 2, and
$\sim C/\sqrt{\log x}$ in Example 3 above).   If there is a
substantial difference between ${\Cal A}(y)$ and $y{\Cal A}(x)/x$
then this already indicates large scale fluctuations
in the distribution of ${\Cal A}$.
\endremark

Corollary 1.3 gives a Roth-type result for general arithmetic
sequences which do not look like the set of all natural numbers.
We will deduce it in Section 2 from the stronger, but more
technical, Theorem 2.4 below.  Clearly Corollary 1.3 applies to
the sequences of primes (with $\alpha=1+o(1)$) and sums of two
squares (with $\alpha=1/2+o(1)$), two results already known.
Surprisingly it applies also to any subset of the primes:

\noindent{\bf Example 4.}  Let ${\Cal A}$ be any subset of the primes.
Then for any fixed $u\ge 1$ and
sufficiently large $x$ there exists $\ell \le x/(\log x)^u$
such that, for some $y\in (x/4,x)$ and some arithmetic progression $a \pmod \ell$
with $(a,\ell)=1$, we have
$$
\Big| {\Cal A}(y;\ell,a)- \frac{1}{\phi(\ell)}\frac{y{\Cal A}(x)}{x}\Big| \gg_u
\frac{{\Cal A}(x)}{\phi(\ell)}.
$$
This implies the first result of section 1a.  
A similar result holds for any subset of the numbers that are sums of two squares.

\noindent{\bf Example 5.} Let ${\Cal A}$ be any subset of those
integers $\leq x$ having no prime factor in the interval $[(\log
x)^{1-\epsilon}, \log x]$. We can apply Corollary 1.3 since
$\alpha\geq \epsilon+o(1)$, and then easily deduce the second
result of section 1a.

Our next result gives an ``uncertainty principle'' implying that
we either have poor distribution in long arithmetic progressions,
or in short intervals.

\proclaim{Corollary 1.4}  Let $\Cal A$, ${\Cal S}$, $h$, $f_q$ and
$\gamma_q$ be as above. Suppose that $0\le h(n)\le 1$ for all $n$.
Suppose that {\rm (1.7)} holds for some $\alpha \ge 60\log \log
\log x/\log \log x$ and set $\eta =\min(\alpha/3,1/100)$. Then for
each $5/\eta^2 \le u\le \eta (\log x)^{\eta/2}$ at least one of
the following two assertions holds:

{\rm (i)} There exists an interval $(v,v+y) \subset (x/4,x)$ with $y\ge (\log x)^u$
such that
$$
\Big|{\Cal A}(v+y) -{\Cal A}(v) - y\frac{{\Cal A}(x)}{x}\Big|
\gg \exp\Big( -\frac{u}{\eta} (1+25\eta) \log (2u/\eta^3) \Big) y
\frac{{\Cal A}(x)}{x}.
$$

{\rm (ii)} There exists $y \in (x/4,x)$ and an arithmetic progression $a \pmod q$
with $(q,{\Cal S})=1$ and  $q\le \exp({2(\log x)^{1-\eta}})$ such that
$$
\Big|{\Cal A}(y;q,a) - \frac{f_q(a)}{q\gamma_q} y\frac{{\Cal A}(x)}{x}
\Big| \gg  \exp\Big( -\frac{u}{\eta} (1+25\eta) \log (2u/\eta^3)
\Big)\frac{{\Cal A}(x)}{\phi(q)}.
$$
\endproclaim

Corollary 1.4 is our general version of Maier's result; it is a
weak form of the more technical Theorem 2.5.  Again condition
(1.7) is invoked to keep away from Example 1. 
Note that we are only able to conclude a dichotomy: either 
there is a large interval $(v,v+y)\subset (x/4,x)$ with 
$y \ge (\log x)^u$ where the density of ${\Cal A}$ is 
altered, or there is an arithmetic progression to 
a very small modulus ($q\le x^\epsilon$) where the distribution 
differs from the expected.  This is unavoidable in general,
and our ``uncertainty principle'' is aptly named, for we can
construct sequences (see \S 6a, Example 6) which are well
distributed in short intervals (and then by Corollary 1.4 such a
sequence will exhibit fluctuations in arithmetic progressions).
In Maier's original result the sequence was easily proved to
be well-distributed in these long arithmetic progressions (and so
exhibited fluctuations in short intervals, by Corollary 1.4).

Our proofs develop Maier's ``matrix method'' of playing off
arithmetic progressions against short intervals or other
arithmetic progressions (see \S 2).  In the earlier work on primes
and sums of two squares,  the problem then reduced to showing
oscillations in certain sifting functions arising from the theory
of the half dimensional (for sums of two squares) and linear (for
primes) sieves.  In our case the problem boils down to proving
oscillations in the mean-value of the more general class of
multiplicative functions satisfying $0\le f(n)\le 1$ for all $n$
(see Theorem 3.1).  Along with our general formalism, this forms
the main new ingredient of our paper and is partly motivated by
our previous work [6] on multiplicative functions and integral
equations.  In \S 7 we present a simple analogue of such
oscillation results for a wide class of integral equations which
has the flavor of a classical ``uncertainty principle'' from
Fourier analysis.

This broader framework has allowed us to improve the uniformity of
the earlier result for primes, and to obtain perhaps best possible
results in this context.

\proclaim{Theorem 1.5}  Let $x$ be large and
suppose $\log x\le y\le \exp(\beta
\sqrt{\log x}/2\sqrt{\log\log x})$, for a certain absolute constant $\beta>0$.
Define
$$
\Delta(x,y) = (\vtheta(x+y)-\vtheta(x)-y)/y,
$$
where $\vtheta(x)=\sum_{p\le x} \log p$.
There exist numbers $x_{\pm}$ in $(x,2x)$ such that
$$
\Delta(x_+, y) \ge y^{-\delta(x,y)} \qquad
\text{and} \qquad \Delta(x_{-},y) \le
-y^{-\delta(x,y)},
$$
where
$$
\delta(x,y)= \frac{1}{\log\log x} \Big( \log \Big(\frac{\log y}{\log\log x}\Big) +
\log\log \Big ( \frac{\log y}{\log\log x}\Big) + O(1) \Big).
$$
\endproclaim

These bounds  are $\gg 1$ if $y=(\log x)^{O(1)}$.
If $y=\exp((\log x)^{\tau})$ for $0<\tau <1/2$ then these bounds are
$\gg y^{-\tau (1+o(1))}$.  Thus we note that the asymptotic, suggested
by probability considerations,
$$
\vtheta(x+y)-\vtheta(x)= y+O(y^{\frac 12+\epsilon}),
$$
fails sometimes for $y\le \exp((\log x)^{\frac 12-\epsilon})$. A.
Hildebrand and Maier [14] had previously shown such a result for
$y\le \exp((\log x)^{\frac 13-\epsilon})$ (more precisely they
obtained a bound $\gg y^{-(1+o(1))\tau/(1-\tau)}$ in the range
$0<\tau <1/3$), and were able to obtain our result assuming the
validity of the Generalized Riemann Hypothesis. We have also been
able to extend the uniformity with which Friedlander and
Granville's result (1.1) holds, obtaining results which previously
Friedlander, Granville, Hildebrand and Maier [4] established
conditionally on the Generalized Riemann Hypothesis.  We will
describe these in \S 5.

This paper is structured as follows:\  In section 2 we describe
the framework in more detail, and show how Maier's method reduces
our problems to exhibiting oscillations in the mean-values of
multiplicative functions.  This is investigated in section 3 which
contains the main new technical results of the paper.  From these
results we quickly obtain in \S 4 our main general results on
irregularities of distribution. In \S 5 we study in detail
irregularities in the distribution of primes.  Our general
framework allows us to substitute a zero-density result of P. X.
Gallagher where previously the Generalized Riemann Hypothesis was
required. In \S 6 we give more examples of sequences covered by
our methods.  Finally in \S 7 we discuss the analogy between
integral equations and mean-values of multiplicative functions,
showing that the oscillation theorems of \S 3 may be viewed as an
``uncertainty principle'' for solutions to integral equations.

\head 2.  The framework \endhead

\noindent  Recall from the introduction that
$a(n)\ge 0$ and that ${\Cal A}(x)=\sum_{n\le x} a(n)$.
Recall that ${\Cal S}$ is a finite set of `bad' primes, and
that $h$ denotes a non-negative multiplicative function
that we shall think of as providing an approximation
$$
{\Cal A}_d(x) :=\sum\Sb n\le x\\ d|n\endSb a(n) \approx
\frac{h(d)}{d} {\Cal A}(x), \tag{2.1}
$$
for each $(d,{\Cal S})=1$.
Roughly speaking, we think of $h(d)/d$ as being the
``probability'' of being divisible by $d$.  The condition 
that $h$ is  multiplicative means that 
 the ``event'' of being divisible by $d_1$ is 
independent of the ``event'' of being divisible by $d_2$, for 
coprime integers $d_1$ and $d_2$.  
We may assume that $h(p^k)<p^k$ for all
prime powers $p^k$ without any significant loss of generality.
As we shall see shortly we may also assume that $h(p^k) \le 1$ without
losing interesting examples.
Let
$$
{\Cal A}(x;q,a) = \sum\Sb n\le x\\ n\equiv a\pmod q\endSb a(n).
$$
We hypothesize that, for $(q,{\Cal S})=1$,
the asymptotics of ${\Cal A}(x;q,a)$ depends only
on the greatest common divisor of $a$ and $q$.  Our aim is to investigate the
limitations of such a model.

First let us describe what (2.1) and our hypothesis predict
for the asymptotics of ${\Cal A}(x;q,a)$.  Writing $(q,a)=m$,
since $|\{ b\pmod q: (b,q)=m\}| = \varphi(q/m)$,
from our hypothesis on ${\Cal A}(x;q,a)$ depending only on $(q,a)$
we would guess that
$$
{\Cal A}(x;q,a)
\approx \frac{1}{\varphi(q/m)} \sum\Sb n\le x \\ (q,n)=m
\endSb a(n)  =\frac{1}{\varphi(q/m)} \sum\Sb n\le x\\ m|n \endSb a(n) \sum\Sb
d| \frac qm\\ d| \frac nm \endSb \mu(d)
= \frac{1}{\varphi(q/m)}\sum_{d |\frac qm } \mu(d){\Cal A}_{dm}(x).
$$
Using now (2.1) we would guess that
$$
{\Cal A}(x;q,a) \approx {\Cal A}(x) \frac{1}{\varphi(q/m)}\sum_{d|\frac qm}
\mu(d) \frac{h(dm)}{dm} =: \frac{f_q(a)}{q\gamma_q} {\Cal A}(x), \tag{2.2}
$$
where
$$
\gamma_q = \prod_{p|q} \Big(\frac{1-h(p)/p}{1-1/p}\Big)^{-1} =
\prod_{p} \Big(1-\frac 1p\Big) \Big(1+\frac{f_q(p)}{p} + \frac{f_q(p^2)}{p^2}
+\ldots \Big) , \tag{2.3}
$$
and $f_q(a)$ is a suitable multiplicative function with
$f_q(a)=f_q((a,q))$ so that it is periodic with period $q$, which
we now define. Evidently $f_q(p^k)=1$ if $p\nmid q$. If $p$
divides $q$, indeed if $p^e$ is the highest power of $p$ dividing
$q$ then
$$
f_q(p^k) := \cases \Big(h(p^k)-\frac{h(p^{k+1})}p\Big)\Big(
1-\frac{h(p)}{p}\Big)^{-1} & \text{if} \ k< e \\
h(p^e)   \Big( 1 - \frac{1}p \Big) \Big(
1-\frac{h(p)}{p}\Big)^{-1} & \text{if} \ k\geq e . \endcases
$$
Note that if $q$ is squarefree and  $h(p)\leq 1$ then
$f_q(p^k)\leq 1$ for all prime powers $p^k$.

We are interested in understanding the limitations to the model (2.2).
We begin with a simple observation that allows us to restrict attention
to the case $0\le h(n)\le 1$ for all $n$.

\proclaim{Proposition 2.1} Suppose that $q\leq x$ is an integer
for which $h(q)>6$. Then either
$$
\Big|{\Cal A}(x;q,0) - \frac{f_q(0)}{q\gamma_q} {\Cal A}(x)\Big|
\ge \frac{1}{2} \frac{f_q(0)}{q\gamma_q} {\Cal A}(x)
$$
or, for every prime $\ell$ in the range $x\ge \ell \ge
3(x+q)/h(q)$ which does not divide $q$, there is an arithmetic
progression $b\pmod{\ell}$ such that
$$
\Big|{\Cal A}(x;\ell, b) - \frac{f_\ell(b)}{\ell \gamma_\ell}
{\Cal A}(x) \Big| \ge \frac 12 \frac{f_\ell(b)}{\ell \gamma_\ell}
{\Cal A}(x).
$$
\endproclaim

The first criterion is equivalent to $|{\Cal A}_q(x) - (h(q)/q)
{\Cal A}(x)| \ge \frac{1}{2} (h(q)/q) {\Cal A}(x)$, since
$f_q(0)/q\gamma_q=f_q(q)/q\gamma_q = h(q)/q$.

\demo{Proof} If the first option fails then
$$
\sum_{n\le x/q} {\Cal A}(x;\ell,nq) \ge \sum_{n\le x/q} a(nq)
={\Cal A}(x;q,0)\ge \frac{1}{2} \frac{f_q(0)}{q\gamma_q} {\Cal
A}(x) =\frac{h(q)}{2q} {\Cal A}(x) .
$$
On the other hand, if prime $\ell \nmid q$ then $f_\ell(nq)=1$ if
$\ell \nmid n$, and $f_\ell(nq)=h(\ell)\gamma_\ell$ if $\ell|n$.
Therefore for any $N$,
$$
\align \sum_{n\leq N} \frac{f_\ell(nq)}{\ell \gamma_\ell} &=
\sum\Sb n\leq N \\ l \nmid n\endSb  \frac{1}{\ell \gamma_\ell} +
\sum\Sb n\leq N
\\ l | n\endSb  \frac{h(\ell)}{\ell} \\
&= \frac{1}{\ell-1} ( [N]-[N/\ell]) - \frac{(\ell-1)}\ell h(\ell)
\{ N/\ell\} \leq \frac{N+1}{\ell} .\\
\endalign
$$
Combining this (taking $N=x/q$) with the display above yields
$$
\sum_{n\le x/q} {\Cal A}(x;\ell,nq) \ge \frac{h(q)}{2q} {\Cal
A}(x)  \ge \frac{3(x+q)}{2q\ell} {\Cal A}(x) \ge
\frac{3}{2}\sum_{n\le x/q} \frac{f_\ell(nq)}{\ell \gamma_\ell}
{\Cal A}(x) ,
$$
which implies the Proposition with $b=nq$ for some $n\le x/q$.
\enddemo

We typically apply this theorem with $h(q) >\log^Ax$ for some
large $A$. This is easily organized if, say, $h(p)\ge 1+\eta$ for
$\ge \eta z/\log z$ primes $p \in (z/2,z)$ where $z\leq \log x$,
and letting $q$ be the product of $[\eta z/\log z]$ of these
primes so that $q=e^{\eta z(1+o(1))}$ and we can select any $\ell$
in the range $x\geq \ell \geq x/\exp( (\eta^2/2) z/\log z)$.

Proposition 2.1 allows us to handle sequences for which $h(p)$ is
significantly larger than $1$ for many primes.  Therefore we will,
from now on, restrict ourselves to the case when $0\le h(n)\le 1$
for all $n$.  Suppose that $(q,{\Cal S})=1$ and define $\Delta_q =
\Delta_q(x)$ by
$$
\Delta_q(x) := \max_{x/4\le y\le x} \ \max_{a \pmod q} \
\Big|{\Cal A}(y;q,a) - \frac{f_q(a)}{q\gamma_q} \frac{y}{x} {\Cal
A}(x)\Big| \ \bigg/ \ \frac{{\Cal A}(x)}{\phi(q)} . \tag{2.4}
$$
In view of (2.2) it seems more natural to consider $|{\Cal
A}(y;q,a) - f_q(a)/(q\gamma_q) {\Cal A}(y)|$ instead of (2.4)
above.  However (2.4) seems to be the most convenient way to
formulate our results, and should be thought of as incorporating a
hypothesis that ${\Cal A}(y)/y$ is very close to ${\Cal A}(x)/x$
when $x/4\le y\le x$. Formally we say that  ${\Cal A}(x)/x$ is
slowly varying: a typical case is when ${\Cal A}(x)/x$ behaves like
a power of $\log x$, a feature seen in the motivating examples of
${\Cal A}$ being the set of primes, or sums of two squares. With
these preliminaries in place we can now formulate our main
principle.

\proclaim{Proposition 2.2}  Let $x$ be large and let ${\Cal A}$,
${\Cal S}$ $h$, $f_q$ and $\Delta_q$ be as above. Let $q \le
\sqrt{x} \le \ell \le x/4$ be positive coprime integers with
$(q,{\Cal S})=(\ell,{\Cal S})=1$. Then
$$
\frac{q}{\phi(q)} \Delta_q(x) +\frac{\ell}{\phi(\ell)}
\Delta_\ell(x) + x^{-\frac 18} \gg \Big| \frac{1}{[x/2\ell]}
\sum_{s\le x/(2\ell)} \frac{f_q(s)}{\gamma_q} - 1\Big|.
$$
\endproclaim

\demo{Proof} Let $R:=[x/(4q)] \ge \sqrt{x}/5$ and $S:= [x/(2\ell)]<\sqrt{x}/2$.
We sum the values of $a(n)$ as $n$ varies over the integers
in the following $R\times S$ ``Maier matrix.''
$$
\boxed{\matrix  (R+1)q+\ell & (R+1)q+2\ell & \cdots & (R+1)q+S\ell \\
 (R+2) q+\ell & (R+2) q+2\ell & \cdots & (R+2) q+S \ell \\
(R+3)q+\ell & \cdot & \cdot & \vdots \\
(R+4) q +\ell & \vdots &
{\boxed{\matrix (r,s)\text{\rm th\ entry}: \\
 (R+r) q + s\ell\endmatrix }} & \vdots \\
\vdots & \vdots & \cdot\ & \vdots  \\
2R q+\ell & \cdots & \cdots  & 2R q+S\ell \endmatrix}
$$
We sum the values of $a(n)$ in two ways: first row by row, and second
column by column.  Note that the $n$ appearing in our ``matrix''
all lie between $x/4$ and $x$.

The $r$-th row contributes
${\Cal A}((R+r)q+\ell S; \ell,(R+r)q)-{\Cal A}((R+r)q;\ell,(R+r)q)$.
Using (2.4), and noting that $f_\ell((R+r)q)=f_\ell(R+r)$ as $(\ell,q)=1$,
this is
$$
\frac{f_\ell(R+r)}{\ell \gamma_\ell} \frac{\ell S}{x}{\Cal A}(x) +
O\Big(\frac{\Delta_\ell}{\phi(\ell)}{\Cal A}(x)\Big).
$$
Summing this over all the rows we see that the sum of $a_n$ with
 $n$ ranging over the Maier matrix above equals
$$
\frac{\ell S}{x} {\Cal A}(x) \sum_{r=R+1}^{2R}
\frac{f_\ell(r)}{\ell \gamma_\ell} + O\Big(
\frac{\Delta_\ell}{\phi(\ell)} {\Cal A}(x) R \Big). \tag{2.5a}
$$

The contribution of column $s$ is ${\Cal A}(2Rq+\ell s;q,\ell s)
-{\Cal A}(Rq+\ell s;q,\ell s)$.  By (2.4), and since $f_q(\ell s)=f_q(s)$ as $(\ell,q)=1$,
we see that this is
$$
\frac{f_q(s)}{q\gamma_q} \frac{Rq}{x} {\Cal A}(x) +
O\Big(\frac{\Delta_q}{\phi(q)} {\Cal A}(x)\Big).
$$
Summing this over all the columns we see that the Maier matrix sum is
$$
\frac{Rq}{x} {\Cal A}(x) \sum_{s=1}^{S} \frac{f_q(s)}{q\gamma_q}
+O\Big(\frac{\Delta_q}{\phi(q)} {\Cal A}(x) S\Big). \tag{2.5b}
$$

Comparing (2.5a) and (2.5b) we deduce that
$$
\frac 1{S\gamma_q} \sum_{s=1}^{S} f_q(s) +
O\Big(\frac{q\Delta_q}{\phi(q)} \Big) = \frac 1{R\gamma_\ell}
\sum_{r=R+1}^{2R} f_\ell(r) + O\Big(\frac{\ell
\Delta_\ell}{\phi(\ell)} \Big). \tag{2.6}
$$

Write $f_\ell(r)=\sum_{d|r} g_\ell(d)$ for a multiplicative function $g_\ell$.
Note that $g_\ell(p^k)=0$ if $p\nmid \ell$.  We also check easily that
$|g_\ell(p^k)|\le (p+1)/(p-1)$ for primes $p |\ell$, and note that $\gamma_\ell =
\sum_{d=1}^{\infty} g_\ell(d)/d$.  Thus
$$
\frac{1}{R\gamma_\ell}\sum_{r=R+1}^{2R} f_\ell(r) =
\frac{1}{R\gamma_\ell} \sum_{d\le 2R} g_\ell(d)\Big(\frac{R}{d}
+O(1) \Big) = 1 + O\Big(\frac{1}{\gamma_\ell}\sum_{d>2R}
\frac{|g_\ell(d)|}{d} + \frac{1}{R\gamma_\ell}\sum_{d\le 2R}
|g_\ell(d)|\Big).
$$
We see easily that the error terms above are bounded by
$$
\ll \frac{1}{R^{\frac 13}\gamma_\ell} \sum_{d=1}^{\infty} \frac{|g_\ell(d)|}{d^{\frac 23}}
\ll \frac{1}{R^{\frac 13}} \prod_{p|\ell} \Big(1+O\Big(\frac{1}{p^{\frac 23}}\Big)\Big)
\ll \frac{1}{R^{\frac 14}},
$$
since $\ell \le x$, and $R\gg \sqrt{x}$.  We conclude that
$$
\frac{1}{R\gamma_\ell}\sum_{r=R+1}^{2R} f_\ell(r) = 1 + O(R^{-\frac 14}).
$$
Combining this with (2.6) we obtain the Proposition.

\enddemo

In Proposition 2.2 we compared the distribution of ${\Cal A}$ in
two arithmetic progressions.  We may also compare the distribution
of ${\Cal A}$ in an arithmetic progression versus the distribution
in short intervals.  Define
$\tilde{\Delta}(y)=\tilde{\Delta}(y,x)$ by
$$
\tilde{\Delta}(y,x) := \max_{(v,v+y)\subset (x/4,x)} \  \Big|{\Cal
A}(v+y)-{\Cal A}(v) - y\frac{{\Cal A}(x)}{x}\Big|\ \bigg/ \
y\frac{{\Cal A}(x)}{x}. \tag{2.7}
$$

\proclaim{Proposition 2.3}  Let $x$ be large and let ${\Cal A}$,
${\Cal S}$, $h$, $f_q$, $\Delta_q$ and $\tilde{\Delta}$ be as
above. Let $q\le \sqrt{x}$ with $(q,{\Cal S})=1$ and let $y\le
x/4$ be positive integers. Then
$$
\frac{q}{\phi(q)} \Delta_q(x) + \tilde{\Delta}(x,y) \gg
\Big|\frac{1}{\gamma_q y} \sum_{s\le y} f_q(s) -1 \Big|.
$$
\endproclaim

\demo{Proof}  The argument is similar to the
proof of Proposition 2.2, starting with an $R\times y$
``Maier matrix'' (again $R=[x/(4q)]$) whose $(r,s)$-th
entry is $(R+r)q+s$.  We omit the details.

\enddemo

We are finally ready to state our main general Theorems
which will be proved in \S 4.

\proclaim{Theorem 2.4}  Let $x$ be large, and in particular
suppose that ${\Cal S}\subset [1,\log \log x]$.
Let $1/100 >\eta \ge 20\log \log \log x/\log \log x$ and suppose that
$(\log x)^{\eta} \le z\le (\log x)/3$ is such that
$$
\sum_{z^{1-\eta} \le p\le z} \frac{1-h(p)}{p} \ge \eta \log
((1-\eta)^{-1}).
$$
Then for all $5/\eta^2 \le u\le \sqrt{z}$
$$
\max\Sb \ell \le x/z^u\\ (\ell, {\Cal S})=1\endSb
 \Delta_\ell \gg \exp\left(-u(1+25\eta)\log(2u/\eta^2)\right).
$$
\endproclaim

Note that $\sum_{z^{1-\eta} \le p\le z}1/p \sim \log
((1-\eta)^{-1})$. There is an analogous result for short
intervals.

\proclaim{Theorem 2.5} Let $x$ be large, and in particular suppose
that ${\Cal S}\subset [1,\log \log x]$.  Let
$1/100 \ge \eta \ge 20\log \log \log x/\log \log x$ and
suppose that $(\log x)^{\eta} \le z\le (\log x)/3$ is such that
$$
\sum_{z^{1-\eta} \le p\le z} \frac{1-h(p)}{p} \ge \eta \log
((1-\eta)^{-1}).
$$
Then for each $5/\eta^2\le u\le \sqrt{z}$ at least one of the following
statements is true:

(i) For $q\le e^{2z}$ which is composed only of primes in
$[z^{1-\eta},z]$ (and so with $(q,{\Cal S})=1$) and such that
$\sum_{p|q} (1-h(p))/p \ge \eta^2$, we have $\Delta_q \gg
\exp(-u(1+25 \eta) \log(2u/\eta^2))$.

(ii) There exists $y\ge z^u$ with $\tilde{\Delta}(y) \gg
\exp(-u(1+25 \eta) \log(2u/\eta^2))$.
\endproclaim

\demo{Deduction of Corollary 1.3} We see readily that there exists
$(\log x)^{\eta} \le z\le (\log x)/3$ satisfying the hypothesis of
Theorem 2.4. Applying Theorem 2.4 (with $u/\eta$ there instead of
$u$) we find that there exists $\ell \le x/z^{u/\eta} \le x/(\log
x)^u$ with $(\ell, {\Cal S})=1$ and $\Delta_\ell \gg
\exp(-(u/\eta)(1+25\eta)\log(2u/\eta^3))$. The corollary follows
easily.

\enddemo

\demo{Deduction of Corollary 1.4}  We may find $(\log x)^{\eta}
\le z\le (\log x)^{1-\eta}$ satisfying the hypothesis of Theorem
2.5. The corollary follows easily by applying Theorem 2.5 with
$u/\eta$ there in place of $u$.
\enddemo

\head 3.  Oscillations in mean-values of multiplicative functions
\endhead

\subhead 3a. Large oscillations \endsubhead

\noindent Throughout this section we shall assume that $z$ is
large, and that $q$ is an integer all of whose prime factors are
$\le z$. Let $f_q(n)$ be a multiplicative function with
$f_q(p^k)=1$ for all $p\nmid q$, and
$0\le f_q(n) \le 1$ for all $n$. Note that $f_q(n)=f_q((n,q))$ is periodic 
$\pmod q$.  Define
$$
F_q(s) = \sum_{n=1}^{\infty} \frac{f_q(n)}{n^s} =
\zeta(s) G_q(s), \qquad \text{where}
\ \ \
G_q(s) = \prod_{p|q} \Big(1-\frac{1}{p^s} \Big) \Big( 1+ \frac{f_q(p)}{p^s}
+\frac{f_q(p^2)}{p^{2s}} +\ldots \Big).
$$
To start with $F_q$ is defined in Re$(s) >1$, but note that the above
furnishes a meromorphic continuation to Re$(s)>0$.
Note also that $\gamma_q=G_q(1)$ in the notation of \S 2.
Define
$$
E(u):= \frac{1}{z^u} \sum_{n\le z^u} (f_q(n) - G_q(1)),
$$
and put for all complex numbers $\xi$
$$
H_j (\xi) := \sum_{p|q}  \frac{1-f_q(p)}{p} p^{\xi/\log z}
\Big(\frac{\log (z/p)}
{\log z}\Big)^{j} \ \text{for each} \ j\geq 0, \ \
\text{and} \ J(\xi):= \sum_{p|q}
\frac{1}{p^2} p^{2\xi/\log z}.
$$
Let $H(\xi):=H_0(\xi)$.

\proclaim{Theorem 3.1}  With notations as above, we
have for $1\le \xi \le \frac 23 \log z$
$$
|E(u)| \le \exp( H(\xi) -\xi u  + 5 J(\xi)).
$$
Let $\frac 23\log z \ge \xi\geq \pi$ and suppose
that $H(\xi) \ge 20 H_2(\xi) + 76 J(\xi) +20$, so that
$$
\tau:= \sqrt{(5H_2(\xi)+19J(\xi)+5)/H(\xi)} \le 1/2.
$$
Then there exist points
$u_{\pm}$ in the interval $[H(\xi)(1-2\tau), H(\xi)(1+2\tau)]$ such that
$$
E(u_+) \ge \frac{1}{20\xi H(\xi)} \exp \{ H(\xi)-\xi u_+ - 5H_2(\xi)-5J(\xi) \},
$$
and
$$
E(u_-) \le - \frac{1}{20\xi H(\xi)} \exp \{ H(\xi) - \xi u_{-} -5H_2(\xi)- 5J(\xi)\}.
$$
\endproclaim

In section 3b (Proposition 3.8) we will show that under certain special circumstances
one can reduce length of the range for $u_\pm$ to 2.

We now record some corollaries of Theorem 3.1.

\proclaim{Corollary 3.2} Let $z^{-\frac 1{10}}\le \eta \le 1/100$ and suppose
that $q$ is composed of primes in $[z^{1-\eta},z]$ and that
$$
\sum_{p|q} \frac{1-f_q(p)}{p} \ge \eta^2.
$$
Then for $\sqrt{z} \ge u \ge 5/\eta^2$ there
exist points $u_\pm \in [u,u(1+22\eta)]$ such that
$$
E(u_+) \ge \exp\Big(-u(1+25\eta) \log \Big(\frac{2u}{\eta^2}\Big)\Big),
\qquad \text{and} \qquad
E(u_-) \le - \exp\Big(-u(1+25\eta) \log \Big(\frac{2u}{\eta^2}\Big)\Big).
$$
\endproclaim
\demo{Proof}  Note that for $1\le \xi \le \frac{11}{20}\log z$
$$
H(\xi) = \sum_{p|q} \frac{1-f_q(p)}{p} p^{\xi/\log z} \ge \eta^2 e^{(1-\eta)\xi},
$$
and that
$$
H_2(\xi) \le \eta^2 H(\xi), \qquad \text{and} \qquad J(\xi) \le
e^{2\xi}{z^{1-\eta}} \le \eta^2 H(\xi),
$$
where the last inequality for $J(\xi)$ is easily checked using our
lower bound for $H(\xi)$ and keeping in mind that $z^{-\frac
1{10}}\le \eta \le 1/100$ and that $\xi \le \frac {11}{20} \log z$. From
these estimates it follows that if $H(\xi)\ge 5/\eta^2$ then
$\tau$ (in Theorem 3.1) is $\le 5\eta$.
Therefore from Theorem 3.1 we conclude that if
$H(\xi)\ge 5/\eta^2$ and $\pi \le \xi \le \frac {11}{20} \log z$ then
there exist points $u_{\pm}$ in
$[H(\xi)(1-10\eta),H(\xi)(1+10\eta)]$ such that
$$
E(u_+) \ge \frac{1}{20\xi H(\xi)} \exp(H(\xi)-\xi u_+-5H_2(\xi)-5J(\xi))
\ge e^{-\xi u_+},
$$
and $E(u_-)\le -e^{-\xi u_-}$.  Renaming $H(\xi)(1-10\eta)=u$
so that $\xi \le \frac{1}{1-\eta}\log (2u/\eta^2)$ we easily
obtain the Corollary.

\enddemo

\proclaim{Corollary 3.3}  Suppose that
$q$ is divisible only by primes between $\sqrt{z}$ and $z$.  Further
suppose $c$ is a positive constant such that for
$1\le \xi \le \frac 23 \log z$ we have $H(\xi) \ge ce^{\xi}/{\xi}$.
Then there is a positive constant $A$ (depending only on $c$) such that
for all $e^A \le u \ll c z^{2/3}/\log z$,
the interval $[u(1-A/\log u), u(1+A/\log
u)]$  contains points $u_{\pm}$ satisfying
$$
E(u_+) \ge \exp\{ -u_+ (\log u_+ + \log\log u_+ +O(1))\},
$$
and
$$
E(u_-) \le - \exp\{ -u_+ (\log u_+ + \log\log u_+ +O(1))\}.
$$
\endproclaim

The implied constants above depend only on $c$. Note that
$\sum_{\sqrt{z}\le p\le z} 1/p^{1-\xi/\log z} \asymp e^{\xi}/\xi$,
by the prime number theorem.  Thus $H(\xi) \ll e^{\xi}/\xi$, and
the criterion $H(\xi) \ge c e^{\xi}/\xi$ in Corollary 3.3 may be
loosely interpreted as saying that, ``typically'', $1-f_q(p) \gg
c$.

 If $H(\xi) \sim \kappa e^{\xi}/{\xi}$ then our bounds take
the shape $\exp\{ -u (\log(u/\kappa) + \log\log u -1 +o(1))\}$.

\demo{Proof of Corollary 3.3} In this situation $J(\xi) \le
\sum_{\sqrt{z}\le p \le z} p^{2\xi/\log z}/p^2 \ll
e^{\xi}/\sqrt{z} +e^{2\xi}/{z}\ll e^{\xi}/\xi^3$.  Further using
the prime number theorem
$$
H_2(\xi) \le \sum_{\sqrt{z} \le p \le z} \frac{p^{\xi/\log z}}{p} \Big( \frac{\log
(z/p)}{\log z}\Big)^2 \ll \int_{\sqrt{z}}^{z} \frac{t^{\xi/\log z}}{t\log t} \Big(
\frac{\log (z/t)}{\log z}\Big)^2 dt \ll \frac{e^{\xi}}{\xi^3}.
$$
The corollary now follows from Theorem 3.1, and renaming $u=H(\xi)$ so that $\xi =
\log u + \log\log u +O(1)$.
\enddemo

\proclaim{Corollary 3.4}  Keep notations as in Theorem 3.1, and suppose $q$ is
divisible only by the primes between $z/2$ and $z$.
Further suppose that $c$ is a positive constant such that for $1\le \xi \le
\frac 23 \log z$ we have $H(\xi) \ge c e^{\xi}/\log z$. Then there is a positive
constant $A$ (depending only on $c$) such that for all $e^A \le u \ll c z^{2/3}/\log
z$ the interval $[u(1-A/\log u), u(1+A/\log u)]$ contains points $u_{\pm}$ satisfying
$$
E(u_+) \ge \frac{1}{\log\log z} \exp\{ - u_+ (\log u_+ + \log\log z+ O(1)) \},
$$
and
$$
E(u_-) \le - \frac{1}{\log\log z} \exp\{ - u_- (\log u_- +
\log\log z+ O(1)) \}.
$$
\endproclaim

As in Corollary 3.3, the implied constants above depend only on $c$.  Also
note that $H(\xi)$ in this case is always $\le \sum_{z/2 \le p\le z}
1/p^{1-\xi/\log z}\asymp e^{\xi}/\log z$.

\demo{Proof of Corollary 3.4} In this case $J(\xi) \ll e^{2\xi}
\sum_{p\geq z/2} 1/p^2  \ll e^{2\xi}/z\log z \ll e^{\xi}/\log^3z$.
Further note that $H_2(\xi) \le (e^{\xi}/\log^2 z)\sum_{z/2 \le
p\le z} 1/p \ll e^{\xi}/\log^3 z$. Taking $u=H(\xi)$ so that $\xi
= \log u +\log\log z+ O(1)$ and thus $H(\xi)\ll e^\xi/\log z$, we
easily deduce Corollary 3.4 from Theorem 3.1.

\enddemo

\proclaim{Corollary 3.5}  Keep the notations of Theorem 3.1, and
suppose (as in Corollary 3.3) that $q$ is
divisible only by primes between $\sqrt{z}$ and $z$ and that
for $1\le \xi \le \frac 23 \log z$ we have $H(\xi) \ge c
e^{\xi}/\xi$.  Let $y=z^u$ with $1\le u\ll cz^{2/3}/\log z$.  There is a positive
constant $B$ (depending only on $c$) such that the interval $[1, z^{u(1+B/\log (u+1))
+B}]$  contains numbers $v_{\pm}$ satisfying
$$
\frac{1}{y} \sum\Sb v_+ \le n \le v_+ +y \endSb (f_q(n)- G_q(1)) \ge
 \exp\{-u (\log (u+1) +\log\log (u+2) +O(1)) \},
$$
and
$$
\frac{1}{y} \sum\Sb v_{-} \le n\le v_{-} + y\endSb (f_q(n)-G_q(1))
 \le  -\exp\{- u (\log (u+1) +\log\log (u+2) +O(1))\}.
$$
\endproclaim

\demo{Proof} Appealing to Corollary 3.3 we see that there is some $w=z^{u_1}$ with
$u_1 \in [e^A + u(1+D/\log(u+1)), e^A + u(1+(D+3A)/\log(u+1))]$ (here $A$ is as in
Corollary 3.3 and $D$ is a suitably large positive constant) such that
$$
\sum\Sb n\le w \endSb (f_q(n)-G_q(1)) \ge
 w \exp( -u_1 (\log u_1 +\log\log u_1 + C_1)) \Big)
$$
for some absolute constant $C_1$.

We now divide the interval $[1,w]$ into subintervals of the form $(w-ky,w-(k-1)y]$ for
integers $1\le k \le [w/y]$, together with one last interval $[1,y_0]$ where $y_0=w-
[w/y]y =y\{w/y\}$. Put $y_0 =z^{u_0}$.  Then, using the first part of Theorem 3.1 to
bound $|E(u_0)|$ (taking there $\xi = \log (u_0+1) + \log\log (u_0+2)$), we get that
$$
\Big |\sum\Sb n\le y_0 \endSb (f_q(n)-G_q(1)) \Big| = y_0 |E(u_0)| \le y_0 \exp( -
(u_0+1) (\log (u_0+1) +\log\log (u_0+2) -C_2))
$$
for some absolute constant $C_2$.

From the last two displayed equations we conclude that
if $D$ is large enough (in terms of $C_1$ and $C_2$) then
$$
\sum\Sb y_0 \le n\le w \endSb (f_q(n)-G_q(1)) \ge w \exp\{ - u (\log (u+1) +\log\log
(u+2) +O(1)) \}
$$
so that the lower bound in the Corollary follows with $v_+ = w- ky$ for some $1\le k
\le [w/y]$.  The proof of the upper bound in the Corollary is similar.
\enddemo

We now embark on the proof of Theorem 3.1. We will write below
$f_q(n)=\sum_{d|n} g_q(d)$ for a multiplicative function $g_q$,
the coefficients of the Dirichlet series $G_q(s)$. Note that
$g_q(p^k)=0$ for $p\nmid q$, and if $p|q$ then $g_q(p^k)=
f_q(p^k)-f_q(p^{k-1})$. Clearly $G_q(1)=\sum_{d=1}^{\infty}
g_q(d)/d$.  Let $[t]$ and $\{t\}$ denote respectively the integer
and fractional part of $t$.  Then
$$
\align E(u) &=  \frac{1}{z^u} \sum_{n\le z^u} \sum_{d|n} g_q(d) - \frac{[z^u]}{z^u}
G_q(1) = \frac{1}{z^u} \sum_{d=1}^{\infty} g_q(d) \Big( \Big[ \frac{z^u}{d}\Big]
- \frac{[z^u]}{d}\Big). \tag{3.1}\\
\endalign
$$
We begin by establishing the upper bound for $|E(u)|$ in Theorem 3.1.

\proclaim{Proposition 3.6}  In the range $1\le \xi \le \frac 23 \log z$ we have
$$
|E(u)| \le \exp( H(\xi) -\xi u  + 5 J(\xi)\Big),
$$
and also
$$
\int_0^{\infty} e^{\xi u} |E(u)| du \le \frac {3}{\xi} \exp( H(\xi) + 5 J(\xi)).
$$
\endproclaim

As will be evident from the proof the condition $\xi \le \frac 23 \log z$ may be
replaced by $\xi \le (1-\epsilon) \log z$.  The constants $3$ and $5$ will have to be
replaced with appropriate constants depending only on $\epsilon$.

\demo{Proof of Proposition 3.6} Since $|[z^u/d] -[z^u]/d|  \le
\min(z^u/d,1)$ we obtain, from (3.1), that
$$
|E(u)| \le \sum_{d\le z^u} \frac{|g_q(d)|}{z^u} + \sum_{d > z^u}
\frac{|g_q(d)|}{d} \le e^{-\xi u} \sum_{d=1}^{\infty}
\frac{|g_q(d)|}{d} d^{\xi/\log z};
$$
and also that
$$
\align \int_0^{\infty} e^{\xi u} |E(u)| du
&\le  \sum_{d=1}^{\infty} |g_q(d)| \Big(  \int_0^{\log d/\log z} \frac{e^{\xi u}}{d}
du + \int_{\log d/\log z}^{\infty} \frac{e^{\xi u}}{z^u}
du \Big)\\
&\le \Big(\frac{1}{\xi}+ \frac{1}{\log z-\xi}\Big)
 \sum_{d=1}^{\infty} \frac{|g_q(d)|}{d} d^{\xi/\log z}.\\
\endalign
$$
Now, as each $|g(p^k)|\leq 1$ and as $2^{1/3}/(2^{1/3}-1) < 5$,
$$
\align
 \sum_{d=1}^{\infty} \frac{|g_q(d)|}{d} d^{\xi/\log z} &\le
\prod_{p|q} \Big( 1+ \frac{1-f_q(p)}{p^{1-\xi/\log z}} +
\sum_{k=2}^{\infty} \frac{1}{p^{k(1-\xi/\log z)}} \Big) \\
&\le \prod_{p|q} \Big(1+ \frac{1-f_q(p)}{p^{1-\xi/\log z}} \Big) \Big(1+ \frac{5}{
p^{2(1-\xi/\log z)}}\Big)
\\
\endalign
$$
since $p\ge 2$ and $\xi \le \frac 23 \log z$. The Proposition follows upon taking
logarithms.
\enddemo

Define
$$
I(s) = \int_0^{\infty} e^{-su} E(u) du.
$$
From Proposition 3.6 it is clear that $I(s)$ converges absolutely
in Re$(s) > -\frac 23 \log z$, and thus defines an analytic
function in this region. Further if Re$(s) > 0$ then
$$
\align I(s) &= \int_0^{\infty} \frac{e^{-su}}{z^u} \sum_{n\le z^u} (f_q(n)  -G_q(1)) =
\sum_{n=1}^{\infty} (f_q(n)-G_q(1)) \int_{\log n/\log z}^{\infty}
 \frac{e^{-su}}{z^u}
du \\
&= \frac{\zeta(1+s/\log z)}{\log z +s} ( G_q(1+s/\log z) - G_q(1) ). \tag{3.2}
\\
\endalign
$$
By analytic continuation this identity continues to hold for all Re$(s) > -\frac 23
\log z$.

\proclaim{Proposition 3.7}  For $1\le \xi \le \frac 23 \log z$
with $z$ sufficiently large, we have
$$
\int_0^{\infty} e^{\xi u} |E(u)| du \ge \frac{\xi}{\xi^2 +\pi^2} \Big( \exp\left\{
H(\xi) - 5H_2(\xi) - 5J(\xi) \right\}- 1\Big).
$$
\endproclaim

\demo{Proof} Taking $s= - (\xi +i \pi)$ in (3.2) we deduce that, since $|G_q(1)| \le
1$,
$$
\int_0^{\infty} e^{\xi u} |E(u)| du \ge |I(s)| \ge  \frac{|\zeta(1+s/\log z)|}{|s+\log
z|} \Big( | G_q(1+s/\log z)| -1\Big).
$$
  From the formula $\zeta(w)= w/(w-1) - w \int_1^{\infty} \{x\} x^{-1-w} dx$,
which is valid for all Re$(w) >0$, we glean that
$$
\frac{|\zeta(1+s/\log z)|}{|s+\log z|} \ge \Big| \text{Re} \Big( \frac{-1}{(\xi+i\pi)}
- \frac 1{\log z} \int_{1}^{\infty} \{x\} x^{-2+\xi/\log z} \cos \Big(\frac{\pi \log
x}{\log z} \Big)dx \Big)\Big|.
$$
For large $z$ and $\xi \le \frac 23 \log z$ we see easily that the integral above is
positive\footnote{In fact, for $z\geq 200$.}, and so we deduce that
$$
|\zeta(1+s/\log z)|/|s+\log z|
\ge \text{Re}(1/(\xi+i\pi)) = \xi/(\xi^2+\pi^2).
$$

Next we give a lower bound for $|G_q(1+s/\log z)|$.  We claim that for $z \ge 101^6$ and
for all primes $p$
$$
\Big|1+\frac{f_q(p)}{p^{1+s/\log z}} + \frac{f_q(p^2)}{p^{2(1+s/\log z)}} +\ldots
\Big| \ge \Big| 1+\frac{f_q(p)}{p^{1+s/\log z}}\Big| \Big( 1-\frac{103}{100}
\frac{1}{p^{2-2\xi/\log z}}\Big). \tag{3.3}
$$
When $p>101^3$ we simply use that the left side of (3.3) exceeds
$|1+f_q(p)/p^{1+s/\log z}| - \sum_{k=2}^{\infty} 1/p^{k(1-\xi/\log
z)}$ and the claim follows. For small $p < 101^3$, set $K=[\log
z/(2\log p)]$ and observe that for $k\leq K$ the numbers
$f_q(p^k)/p^{k(1+s/\log z)}$ all have argument in the range
$[0,\pi/2]$.  Hence the left side of (3.3) exceeds, writing
$q=1/p^{1-\xi/\log z}$,
$$
\Big|1+\frac{f_q(p)}{p^{1+s/\log z}}\Big| - \sum_{k>K} \frac{1}{q^k} \geq
\Big|1+\frac{f_q(p)}{p^{1+s/\log z}}\Big| \left( 1- \frac 1 {q^{K-1}(q-1)^2} \right),
$$
which implies (3.3) for $z\geq 101^6$.

Observe that if $|w|\le 2^{-1/3}$ then
$$
\log |1+w| = \text{Re }
(w - \sum_{n=2}^{\infty} (-1)^n w^n/n) \ge \text{Re }(w) - 5|w|^2/4.
$$
From this observation and our claim (3.3) we deduce easily that
$$
\log \Big| 1-\frac{1}{p^{1+s/\log z}}\Big| \Big| 1+\frac{f_q(p)}{p^{1+s/\log z}}
+\ldots \Big| \ge \text{Re } \left(\frac{f_q(p)-1}{p^{1+s/\log z}}\right) -
\frac{5}{p^{2(1-\xi/\log z)}}.
$$
It follows that
$$
\align \log |G_q(1+s/\log z)| &\ge - \text{Re } \sum_{p|q}
\left(\frac{1-f_q(p)}{p}\right) p^{-s/\log z} -5 J(\xi)
\\
&= H(\xi) + \sum_{p|q} \left(\frac{1-f_q(p)}{p}\right) p^{\xi/\log z} \Big(-1 - \cos
\Big(\frac{\pi\log p}{\log z}\Big) \Big)
-5J(\xi).\\
\endalign
$$
Since $-1-\cos (\pi \log p/\log z) \ge - (\pi^2/2) (\log (z/p)/\log z)^2$, we deduce
the Proposition.

\enddemo


We are now ready to prove Theorem 3.1.

\demo{Proof of Theorem 3.1} The first part of the result was proved in Proposition
3.6. Now, let $I^+$ (and $I^-$) denote the set of values $u$ with $E(u) \ge 0$
(respectively $E(u) <0$). Taking $s=-\xi$ in (3.2) we deduce that for $\xi \le \frac
23\log z$
$$
\Big| \int_0^{\infty} e^{\xi u} E(u) du \Big| \le 2 \Big| \frac{\zeta(1-\xi/\log
z)|}{\log z -\xi}   \Big| \le \frac{6}{\xi},
$$
since $0\leq G_q(1-\xi/\log z), G_q(1)\le 1$, and since (using
$\zeta(w)/w = 1/(w-1)-\int_1^{\infty} \{x\}x^{-1-w} dw$)
$$
\Big|\frac{\zeta(1-\xi/\log z)}{(\log z-\xi)}\Big| = \Big|-\frac{1}{\xi}
-\frac{1}{\log z} \int_1^{\infty} \{x\} x^{-2+\xi/\log z} dx \Big| \le \frac{1}{\xi} +
\frac{1}{\log z-\xi} \le \frac{3}{\xi}.
$$
Combining this with Proposition 3.7 we deduce easily that
$$
\align \int_{I^{\pm}} e^{\xi u} |E(u)| du &\ge \frac{\xi}{2(\xi^2+\pi^2)}
\Big(\exp\{H(\xi)- 5H_2(\xi)-5J(\xi)\} -1 \Big) - \frac{3}{\xi}\\
&\ge \frac{1}{5\xi} \exp\{ H(\xi)-5H_2(\xi)-5J(\xi)\}.
\tag{3.4}\\
\endalign
$$

Put $u_1=H(\xi)(1+2\tau)$.  Then by Proposition 3.6 we get that
$$
\align \int_{u_1}^{\infty} e^{\xi u} |E(u)| du &\le e^{-\tau u_1} \int_0^{\infty}
e^{(\xi+\tau)u} |E(u)| du
\le \frac{3}{\xi} \exp \{ H(\xi+\tau)- \tau u_1 + 5J(\xi+\tau)\}.\\
\endalign
$$
Now $H(\xi+\tau) \le e^{\tau} H(\xi) \le (1+\tau +\tau^2) H(\xi)$, and $J(\xi+\tau)
\le e^{2\tau} J(\xi) \le 2.8 J(\xi)$.   Hence $H(\xi+\tau)+5J(\xi+\tau)-\tau u_1 \le
H(\xi) - 5H_2(\xi) -5J(\xi) -5$, and so we conclude that
$$
\int_{u_1}^{\infty} e^{\xi u} |E(u)| du \le \frac{1}{20\xi} \exp\{ H(\xi) -
5H_2(\xi)-5J(\xi) \}. \tag{3.5}
$$

Similarly note that, with $u_0= H(\xi)(1-2\tau)$,
$$
\int_0^{u_0} e^{\xi u} |E(u)| du \le e^{\tau u_0} \int_0^{\infty} e^{(\xi-\tau)u }
|E(u)| du \le \frac{3}{\xi-\tau}  \exp\{ H(\xi -\tau) + \tau u_0 +  5J(\xi-\tau) \}.
$$
Now $J(\xi-\tau) \le J(\xi)$, and
$$
\align H(\xi-\tau) &=\sum_{p|q} \frac{1-f_q(p)}{p} p^{\xi/\log z} p^{-\tau/\log z} \le
\sum_{p|q} \frac{1-f_q(p)}{p} p^{\xi/\log z} \Big( 1-\tau \frac{\log p}{\log z} +
\frac{\tau^2}{2} \Big)
\\
&= H(\xi) (1-\tau +\tau^2/2) + \tau H_1(\xi) \le H(\xi)(1-\tau +\tau^2/2) + \tau
\sqrt{H(\xi) H_2(\xi)},
\\
\endalign
$$
since $H_1(\xi )\le \sqrt{H(\xi)H_2(\xi)}$ by Cauchy's inequality.
  From these observations and our
definition of $\tau$ it follows that
 $H(\xi-\tau) + 5J(\xi-\tau) + \tau u_0 \le H(\xi) - 5H_2(\xi)
-5J(\xi) -5$, and so
$$
\int_0^{u_0} e^{\xi u} |E(u)| du \le \frac{1}{20 \xi} \exp\{ H(\xi) - 5 H_2(\xi) -
5J(\xi)\}. \tag{3.6}
$$

Combining (3.4), (3.5), and (3.6) we deduce that
$$
\int_{I^{\pm} \cap [u_0,u_1]} e^{\xi u} |E(u)| du \ge \frac{1}{10\xi} \exp\{ H(\xi) -
5 H_2(\xi) -5J(\xi) \}.
$$
Now $u_1-u_0 \leq 4\tau H(\xi) \leq 2 H(\xi)$, so the Theorem follows.

\enddemo

\subhead 3b. Localization of sign changes of $E$ \endsubhead

\noindent We saw in Corollary 3.3 that (in typical situations) 
$E$ changes sign in
intervals of the form $[u(1-A/\log u),u(1+A/\log u)]$.  We
consider now the problem of providing a better localization of the
sign changes of $E$ for small values of $u$. Our main result of this 
section is the following:

\proclaim{Proposition 3.8} With notation as above, suppose that
 $\max_{x\geq u} |E(x)| \gg 1/(G_q(1)\log z))$ for some
$u\geq 6$. Then there exist points $u_+, u_- \in [u-1,u+1]$ such
that $E(u_+), -E(u_-) \geq \max_{x\geq u} |E(x)|$.
\endproclaim

Proposition 3.8 (which may be easily deduced from the Lemmas 
of this section) can be used to reduce the size of the interval in
Theorem 3.1.  In Corollary 3.3 this is simple to state: For $e^A
\le u\le \log \log z/(2\log \log \log z)$ the interval $[u-1,u+1]$
contains points $u_{\pm}$ satisfying the conclusions of Corollary
3.3.

\proclaim{Lemma 3.9} Uniformly for $u>0$ we have
$$
uE(u) +\int_u^{\infty} E(t) dt + \frac{1}{\log z} \sum_{p\le z}
\frac{1-f_q(p)}{p} \log p \ E \Big(u-\frac{\log p}{\log z}\Big) =
O\Big(\frac{1}{G_q(1)\log z}\Big).
$$
\endproclaim

\demo{Proof} Let $E_1(u):= \sum_{d>z^{u}} g_q(d)/d$; and note that
$|E_1(u)| \leq \sum_d |g_q(d)|/d \leq 1/G_q(1)$. By a result of R.R.
Hall (see [13], or (4.1) of [10]) we see that
$$
\frac{1}{z^u} \sum_{d\le z^u} |g_q(d)| \ll \frac{1}{u\log z}
\sum_{d} \frac{|g_q(d)|}{d} \ll \frac{1}{G_q(1)u\log z}.
$$
Therefore, from (3.1) we deduce that
$$
E(u) = - (1+O(z^{-u})) E_1(u) +O\Big(\frac{1}{z^u} \sum_{d\le z^u}
|g_q(d)|\Big)  = -E_1(u) +O \Big(\frac{1}{G_q(1)u\log z}\Big).
\tag{3.7}
$$
Manipulation of $E_1(u)$ yields our identity. The starting point
is the observation that
$$
uE_1(u) + \int_u^{\infty} E_1(t)dt = uE_1(u) + \sum_{d>z^u}
\frac{g_q(d)}{d} \Big(\frac{\log d}{\log z} -u\Big) = \sum_{d>z^u}
\frac{g_q(d)}{d} \frac{\log d}{\log z}. \tag{3.8}
$$
We approximate the left side of (3.8) as follows:
$$\align
\Big|(uE_1(u) + & \int_u^{\infty} E_1(t)dt)+(uE(u) +\int_u^{\infty}
E(t) dt)\Big| \\
&\leq u|E_1(u)+E(u)| + \int_u^{\infty} |E_1(t)+E(t)|dt
\\
& \ll \frac{1}{G_q(1)\log z} + \int_u^{\infty} z^{-t} \Big(
\frac{1}{G_q(1)} + \sum_{d\le z^u} |g_q(d)| \Big)      dt \\
& \ll \frac{1}{G_q(1)\log z} + \frac{1}{\log z} \Big( \sum_{d }
\frac{|g_q(d)|}{d} \Big) \ll \frac{1}{G_q(1)\log z} .
\endalign
$$
Now $\log d =\sum_{m|d }\Lambda(m)$ so that the right side of
(3.8) equals
$$
\frac{1}{\log z} \sum_{m} \Lambda(m) \sum\Sb d>z^u \\ m|d \endSb \frac{g_q(d)}{d}.
$$
The sum over $m$'s above can be restricted to prime powers $p^k$
for $p\le z$ (else $g_q(d)=0$).  Further the contribution of prime
powers $p^k$ with $k\ge 2$ is bounded by $\ll 1/(G_q(1)\log z)$.
Finally for $m=p \le z$ we see that
$$
\sum\Sb d>z^u \\ m|d \endSb \frac{g_q(d)}{d} = \frac{g_q(p)}{p} \sum_{d>z^u/p} \frac{g_q(d)}{d}
+ O\Big(\frac{1}{p^2G_q(1)}\Big)
= -\frac{1-f_q(p)}{p} E_1\Big(u-\frac{\log p}{\log z}\Big)
+ O\Big(\frac{1}{p^2G_q(1)}\Big).
$$
Therefore, by (3.7), this, taken with the estimate for the left
side of (3.8), yields the result.
\enddemo

We call a point $w$ {\sl special} if $|E(w)| = \max_{x\ge w}
|E(x)|$. Since $E(x)\to 0$ as $x\to \infty$ we see that there are
arbitrarily large special points.

\proclaim{Lemma 3.10}  Given $u\ge 2$ either $E(x) =O(1/(G_q(1)\log
z))$ for all $x\ge u$ or there is a special point in $[u,u+1]$.
\endproclaim

\demo{Proof}  Let $w$ denote the infimum of the set of special
points at least as large as $u$, and assume $w> u+1$ (that is,
there is no special point in $[u,u+1]$). Note that $|E(w)| \ge
|E(x)| +O(z^{-u})$ for any $x\ge u$.  If $E(x)$ maintains the same
sign for all $x\ge w$ set $v=\infty$; otherwise let $v$ denote the
infimum of those points $x\geq w$ for which $E(x)$ has the
opposite sign to $E(w)$.  Note that $E(v) =O(1/z^v)$. Taking Lemma
3.9 with $u=w$ and $u=v$ and subtracting we find that
$$
\align wE(w) + \int_w^v E(t)dt &+ \frac{1}{\log z} \sum_{p\le z}
\frac{1-f_q(p)}{p} \log p \  \Big( E\Big(w-\frac{\log p}{\log
z}\Big) - E\Big(v-\frac{\log p}{\log z}\Big)\Big)
\\
&= O\Big(\frac{1}{G_q(1)\log z}\Big).  \tag{3.9}
\endalign
$$
Since $E(t)$ maintains the same sign throughout $[w,v]$ we have
that
$$
\Big| wE(w) + \int_w^v E(t)dt\Big| \ge w|E(w)|,
$$
while on the other hand
$$
\align \Big|\frac{1}{\log z} \sum_{ p\le z} &\frac{1-f_q(p)}{p}
\log p \  \Big( E\Big(w-\frac{\log p}{\log z}\Big) -
E\Big(v-\frac{\log p}{\log z}\Big)\Big) \Big|
\\
&\le \frac{2}{\log z} \sum_{p\le z}  \frac{\log p}{p} \max_{\xi\ge
w-1} |E(\xi)| \le (2+o(1)) \max_{\xi\ge u}  |E(\xi)|
\le (2+o(1)) |E(w)|,\\
\endalign
$$
since $w$ was assumed to be larger than $u+1$. We conclude that
$(u-1+o(1)) |E(w)| = O(1/(G_q(1)\log z))$ which establishes the
Lemma.
\enddemo

\proclaim{Lemma 3.11} If $u\ge 2$ and $E(x)$ maintains the same
sign throughout $[u,u+2]$ then $E(x)=O(1/(G_q(1)\log z))$ for all
$x\ge u+1$.
\endproclaim

\demo{Proof} Suppose not.  Let $w$ denote the infimum of the set
of all special points at least as large as $u+1$.  By Lemma 3.10 we
know that $w\le u+2$. Let $v$ denote the infimum of points $x\ge
u+2$ such that $E(x)$ has the opposite sign to $E(w)$; if no such
point exists set $v=\infty$.  Now $E$ maintains the same sign in
$[w-1,v]$ (since it is a subinterval of $[u,v]$) and so
$$
\Big|wE(w) + \int_{w}^{v} E(t) dt + \frac{1}{\log z}\sum_{p\le z}
\frac{1-f_q(p)}{p} \log p \ E\Big(w-\frac{\log p}{\log
z}\Big)\Big| \ge w|E(w)|;
$$
on the other hand
$$
\Big|\frac{1}{\log z}\sum_{p\le z} \frac{1-f_q(p)}{p} \log p \
E\Big(v-\frac{\log p}{\log z}\Big) \Big| \le |E(w)| (1+o(1)),
$$
since $|E(v-\log p/\log z)|\le \max_{x\ge v-1} |E(x)| \le |E(w)|$.
Therefore, by (3.9), we deduce that $(w-1+o(1)) |E(w)| \leq
O(1/(G_q(1)\log z))$ which proves the Lemma.
\enddemo

\proclaim{Proposition 3.12} Fix $\epsilon>0$. Given a special point
$w>4+\epsilon$ either there exists a point $\xi\in [w-1,w]$ for
which $E(\xi)$ and $E(w)$ have opposite signs and $|E(\xi)| \ge
(w-4-\epsilon) |E(w)|$, or $E(w)=O(1/(G_q(1)\log z))$.
\endproclaim
\demo{Proof} Select $v \in [w+1,w+3]$ with $E(v) = O(z^{-v})$ (if
$v$ does not exist then the Proposition follows from Lemma 3.11).
Choose $\delta=\pm 1$ so that $\delta E(w)>0$. Suppose that
$\delta E(x)> -(2w-v-1-\epsilon) \delta E(w)$ for all $x\in
[w-1,w]$. Then $\delta$ times the right side of (3.9) is
$$
\geq \delta E(w) ( w - \int_w^v 1dt - \frac{1}{\log z} \sum_{p\le
z} \frac{\log p}{p}  \ (2w-v-1-\epsilon+1) ) \geq \epsilon
|E(w)|/2
$$
since each $v-\log p/\log z\geq w$ so that $|E(v-\log p/\log
z)|\le |E(w)|$. The result follows from (3.9) since $2w-v-1\geq
w-4$.
\enddemo

Given $u>3+\epsilon$ we can take $w$ to be a special point in
$[u+1,u+2]$ and then $\xi\in [u,u+2]$.

\head 4. Proof of Theorems 2.4 and 2.5 \endhead

\noindent We shall only provide the proof of Theorem 2.4, the
proof of Theorem 2.5 being entirely similar. We take $q$ to be the
product of all the primes in $(z^{1-\eta},z)$ so that $q\le
\prod_{p\le z}p \leq x^{1/3+o(1)}$. Since ${\Cal S}\subset [1,\log
\log x]$ we also know that $(q,{\Cal S})=1$.  Since $q$ is
squarefree we check that $f_q(n)$ (as defined in Section 2)
satisfies $f_q(p^k)=1$ if $p\nmid q$ and $f_q(p^k) =
h(p)(p-1)/(p-h(p))$ if $p|q$ and $k\ge 1$ (note that $f_q(p^k)\leq
h(p)\leq 1$ in this case) . Thus $0\le f_q(n) \le 1$ for all $n$,
and we may apply the results of Section 3, in particular Corollary
3.2. From Corollary 3.2 we obtain that for $\sqrt{z}\ge u\ge
5/\eta^2$ there exists $\lambda \in [u,u(1+22\eta)]$ such that
$$
\Big| \frac{1}{z^{\lambda}} \sum_{n\le z^{\lambda}} f_q(n) - G_q(1) \Big|
 \ge \exp\Big( - u(1+25 \eta) \log \Big(\frac{2u}{\eta^2}
\Big)\Big). \tag{4.1}
$$

We now use Proposition 2.2.  We may easily find $\ell \in (x/(2(z^{\lambda}+1)),
x/(2z^{\lambda}))$
such that $\ell$ is not divisible by any prime below $(\log x)/3$.  Such an
$\ell$  is coprime to $q$, satisfies $(\ell,{\Cal S})=1$,
and further $\phi(\ell)/\ell \gg 1$.  Proposition 2.2
therefore yields, by (4.1),
$$
\Delta_q+ \Delta_\ell \gg \exp\Big( - u(1+25 \eta) \log
\Big(\frac{2u}{\eta^2} \Big)\Big),
$$
proving Theorem 2.4 (in the statement of which we might take
$\ell=q$).

\head 5. Limitations to the equidistribution of primes \endhead

\noindent In this section we will exploit the flexibility afforded
by our general oscillation results to obtain refinements to
previous results on the limitations to the equidistribution of
primes.  We stated in Theorem 1.5 (see Introduction) the result for primes
in short intervals and we now state the analogous results
for primes in arithmetic progressions.

\proclaim{Theorem 5.1}  Let $\ell$ be large and suppose that
$\ell$ has fewer than $(\log \ell)^{1-\epsilon}$ prime divisors
below $\log \ell$. Suppose that $(\log \ell)^{1+\epsilon} \le y\le
\exp(\beta \sqrt{\log \ell} /\sqrt{2\log\log \ell})$ for a certain
absolute constant $\beta>0$, and put $x=\ell y$.  Define for
integers $a$ coprime to $\ell$
$$
\Delta(x;\ell,a) = \left( \vtheta(x;\ell,a) - \frac{x}{\vphi(\ell)} \right) \bigg/
\frac{x}{\vphi(\ell)} .
$$
There exist numbers $x_{\pm}$ in the interval $(x,x y^{D/\log (\log y/\log\log
\ell)})$, and integers $a_{\pm}$ coprime to $\ell$ such that
$$
\Delta(x_{+}; \ell, a_+) \ge y^{-\delta(\ell,y)}, \qquad \text{and} \qquad
\Delta(x_-;\ell,a_-) \le -y^{-\delta(\ell,y)}.
$$
Here $D$ is an absolute positive constant which depends  only on $\epsilon$,
and $\delta(\cdot,\cdot)$ is as in Theorem 1.5.
\endproclaim

These bounds are $\gg 1$ if $y=(\log \ell)^{O(1)}$, and $\gg
y^{-\tau (1+o(1))}$ if $y=\exp((\log \ell)^{\tau})$ with $0\le
\tau <1/2$. The corresponding result in [4], Theorem A1, gives the
weaker bound $y^{-(1+o(1))\tau/(1-\tau)}$ (though our bound is
obtained there assuming the Generalized Riemann Hypothesis). Our
constraint on the small primes dividing $\ell$ is less restrictive
than the corresponding condition there, though our localization of
the $x_{\pm}$ values is worse (in [4] the $x_{\pm}$ values are
localized in intervals $(x/2,2x)$).

Theorem 5.1 omits a thin set of moduli $\ell$ having very many small prime factors.  We
next give a weaker variant which includes all moduli $\ell$.

\proclaim{Theorem 5.2}  Let $\ell$ be large and suppose that
$(\log \ell)^{1+\epsilon} \le y \le \exp(\beta \sqrt{\log
\ell}/\sqrt{ \log\log \ell} )$ for a certain absolute constant
$\beta>0$, and put $x=y\ell$. There exist numbers $x_{\pm}$ in the
interval $(x,x y^{D/\log (\log y/\log\log \ell)})$, and integers
$a_{\pm}$ coprime to $\ell$ such that
$$
\Delta(x_{+}; \ell, a_+) \ge \frac{y^{-\delta_1(\ell,y)}}{\log\log\log \ell} \qquad
\text{and} \qquad \Delta(x_-;\ell,a_-) \le -\frac{y^{-\delta_1(\ell,y)}}{\log\log\log
\ell},
$$
where $\delta_1(x,y)= (\log\log y + O(1))/(\log\log x)$. Here $D$ is an absolute
positive constant which depends  only on $\epsilon$.
\endproclaim

Theorem 5.2 should be compared with Theorem A2 of [4]. Our bound is $\gg y^{-\tau
(1+o(1))}$ if $y=\exp((\log \ell)^{\tau})$ with $0<\tau <1/2$. The corresponding
result in [4], Theorem A2, gives a weaker bound $y^{-(3+o(1))\tau/(1+\tau)}$, though
our localization of the $x_{\pm}$ values is again much worse.

To prove Theorems 1.5, 5.1 and 5.2 we require knowledge of the distribution of
primes in certain arithmetic progressions.  We begin by describing such a result,
which will be deduced as a consequence of a theorem of Gallagher [5].

For $1\le j \le J:=[\log z/(2\log 2)]$ and consider the dyadic intervals
$I_j =(z/2^j,z/2^{j-1}]$.  Let $P_j$ denote a subset of the primes in $I_j$, and let
$\pi_j$ denote the cardinality of $P_j$.  We let ${\Cal Q}$ denote the set of integers
$q$ with the following property:  $q = \prod_{j=1}^{J} q_j$ and each $q_j$ is the
product of exactly $[\pi_j/2]$ distinct primes in $P_j$.  It is clear that all the
elements of ${\Cal Q}$ are squarefree and lie below $Q:= z^{\sum_{j} \pi_j/2}$, and
that $|{\Cal Q}| = \prod_{j=1}^{J} \binom{\pi_j}{[\pi_j/2]}$.

There is a constant $c_1$ such that at most one primitive $L$-function with modulus
$q$ between $\sqrt{T}$ and $T$ has a zero in the region $\sigma > 1-c_1/\log q$, and
$|t| \le T$.  Further if this exceptional Siegel zero exists then it is real, simple
and unique (see chapter 14 in [2]).  We call the modulus of such an exceptional
character a Siegel modulus. Below $Q$ there are $\ll \log \log Q$ Siegel moduli.
Denote these by $\nu_1, \nu_2, \ldots, \nu_\ell$, and for each select a prime divisor
$v_1,\dots, v_\ell$. Assume none of $v_1,\dots v_\ell$ belongs to $\cup_{j=1}^{J} P_j$,
which  guarantees that there are no Siegel zeros for any of the moduli $d$ where $d$
is a divisor of some $q\in {\Cal Q}$.

\proclaim{Proposition 5.3}  Suppose that $\exp(\sqrt{\log x}) \le Q \le x^b$ where $b$
is a positive absolute constant, and let $x\exp(-\sqrt{\log x}) \le h\le x$. Then
$$
\frac{1}{|{\Cal Q}|} \sum_{q\in {\Cal Q}} \max\Sb (a,q)=1 \endSb \left|  \frac
{\vtheta(x+h;q,a)-\vtheta(x;q,a) - h/\vphi(q)} {h/\vphi(q)} \right|\ll
\exp\Big(-\alpha \frac{\sqrt{\log x}} {\sqrt{\log z}}\Big),
$$
where $\alpha$ is a positive absolute constant.
\endproclaim

\demo{Proof} If $(a,q)=1$ then using the orthogonality of characters
$$
\align \vtheta(x+h;q,a)-\vtheta(x;q,a) &=\frac{1}{\vphi(q)} \sum_{\chi \pmod q}
\overline{\chi(a)}
\sum_{x}^{x+h} \chi(p) \log p \\
&= \frac{1}{\vphi(q)} \sum\Sb x\le p\le x+h\\
(p,q)=1\endSb \log p + O\Big( \frac{1}{\vphi(q)} \sum_{\chi \neq \chi_0} \Big|
\sum_{x}^{x+h} \chi(p)\log p\Big| \Big).
\\
\endalign
$$
By the prime number theorem the first term above is $h/\vphi(q)
\{1+O(\exp(-c\sqrt{\log x}))\}$, which has an acceptable error term.  We now focus on
estimating the second term on average.

Below the superscript $*$ will indicate a restriction to primitive characters.
Observe that
$$
\align \sum_{\chi \neq \chi_0} \Big| \sum_{x}^{x+h} \chi(p)\log p\Big| &= \sum\Sb d|q
\\ d>1 \endSb\,\, \, \,\,\,
 \sumstar_{\chi \pmod d}
\Big( \Big| \sum_{x}^{x+h} \chi(p)\log p\Big|  +O\Big( \sum_{p | \frac qd} \log
p\Big)\Big)
\\
&=  \sum\Sb d|q \\ d>1 \endSb \,\, \, \,\,\, \sumstar_{\chi \pmod d} \Big|
\sum_{x}^{x+h} \chi(p)\log p\Big| +O(d(q)\log q).
\\
\endalign
$$
Thus
$$
\sum_{q\in {\Cal Q}} \sum_{\chi \neq \chi_0} \Big| \sum_{x}^{x+h} \chi(p)\log p\Big| =
\sum_{1< d\le Q} \Big( \sum\Sb q\in {\Cal Q}\\ d|q\endSb 1\Big) \sumstar_{\chi\pmod d}
\Big| \sum_{x}^{x+h} \chi(p)\log p\Big| + O(|{\Cal Q}| x^{\epsilon}). \tag{5.1}
$$

Observe that if $d$ is to have any multiples in ${\Cal Q}$ then we must have
$d=\prod_{j=1}^{J} d_j$ with each of the $d_j$ being composed of at most $[\pi_j/2]$
distinct primes from $P_j$.  Therefore
$$
\sum\Sb q\in {\Cal Q}\\ d|q\endSb 1\leq \prod_{j=1}^{J}
\binom{\pi_j-\omega(d_j)}{[\pi_j/2]- \omega(d_j)} \le \prod_{j=1}^{J}
2^{-\omega(d_j)} \binom{\pi_j}{[\pi_j/2]} = 2^{-\omega(d)}|{\Cal Q}| \le
2^{-\frac{\log d}{\log z}} |{\Cal Q}|,
$$
where $\omega(n)$ denotes the number of prime factors of $n$, and the final estimate
follows since all prime divisors of $d$ are below $z$ and so $\omega(d) \ge \log
d/\log z$.  From these remarks we see that the right side of (5.1) is
$$
\le |{\Cal Q}| \sumflat_{1<d\le Q} 2^{-\frac{\log d}{\log z}} \sumstar_{\chi \pmod d}
\Big| \sum_{x}^{x+h} \chi(p)\log p\Big| + O( |{\Cal Q}| x^{\epsilon}), \tag{5.2}
$$
where the $\flat$ on the sum over $d$ indicates that the sum is restricted to $d$ as
above; note that such $d$ have no Siegel zeros.

Define $J_0:=[1,\exp(\sqrt{\log x})]$, and  $J_k:=(\exp(\sqrt{\log x})
z^{k-1},\exp(\sqrt{\log x}) z^k]$ for $1\leq k\leq [(\log Q-\sqrt{\log x})/\log z]$.
Theorem 7 of Gallagher [5] implies that the contribution of terms $d\in J_k$ (for $k\ge 1$) is (for
a positive absolute constant $\alpha$)
$$
\le  |{\Cal Q}|
 h \exp\Big(-\alpha^2 \frac{\log x}{\sqrt{\log x} + (k+1)\log z} -
\frac{\sqrt {\log x} +k\log z}{2\log z}\Big) \ll h|{\Cal Q}|
\exp\Big(-\alpha\frac{\sqrt{2\log x}}{\sqrt{\log z}}\Big),
$$
and also that the contribution of the terms
$d\in J_0$ is $\ll |{\Cal Q}| h \exp(-\alpha \sqrt{\log x})$.
Summing over $k$ we deduce that the quantity in (5.2)
is $\ll |{\Cal Q}| h \exp(-\alpha \sqrt{\log x}/\sqrt{\log z})$
which proves the Proposition.

\enddemo

\proclaim{Corollary 5.4}  Suppose that $\exp(\sqrt{\log x}) \le Q \le x^b$ where $b$
is a positive absolute constant. Select $z$ to be the largest integer such that
$\prod_{p\leq z} p \leq Q^{2-o(1)}$. Then there exists an integer $q\in
[Q^{1-c/\log\log Q}, Q]$, whose prime factors all lie in $[\sqrt{z},z]$ with
$\sum_{p|q} 1/p^{1-\xi/\log z} \ge e^{\xi}/3{\xi}$ for $1\le \xi \le \frac 23 \log z$,
such that
$$
\Big| \vtheta(2x;q,a)-\vtheta(x;q,a) - \frac{x}{\vphi(q)}\Big| \ll \frac{x}{\vphi(q)}
\exp\Big(-\beta \frac{\sqrt{\log x}} {\sqrt{\log\log x}}\Big), \tag{5.3}
$$
for all $(a,q)=1$, where $\beta$ is a positive absolute constant.
\endproclaim

\demo{Proof} When $z$ is chosen as above we have $z\sim 2\log Q$ so that $\sqrt{\log
x} \ll z \ll \log x$. If $z/2^j<\sqrt{z}$ then let $P_j=\emptyset$. If $z\geq
z/2^j\geq \sqrt{z}$ then let $P_j$ be the set of primes in $I_j$, omitting the prime
divisors $v_1,\dots v_\ell$ of Siegel moduli: in these cases $\pi_j = \pi(z/2^{j-1})
-\pi(z/2^{j}) + O(\log\log x)$. By Proposition 5.3 we can select $q$ satisfying (5.3),
with $\sim z/(2^{j+1} \log(z/2^j))$ prime factors in each $I_j$ where $j\leq [\log
z/\log 4]$, so that $\sum_{p|q} 1/p^{1-\xi/\log z}  \ge e^{\xi}/3{\xi}$ for $1\le \xi
\le \frac 23 \log z$.
\enddemo

\proclaim{Corollary 5.5}  Fix $1/2\geq \epsilon>0$. Suppose we are given a large
integer $\ell$ with fewer than $(\log \ell)^{1-\epsilon}$ prime divisors below $\log
\ell$. Also we are given $\exp(3(\log \ell)^{1-\epsilon}) \le Q \le \ell^b$, and $J$ a
positive integer $\leq \exp(\sqrt{\log \ell})$. Select $z$ to be the largest integer
such that $\prod_{p\leq z} p \leq Q^{2-o(1)}$. Then there exists an integer $q\in
[Q^{1-c/\log\log Q}, Q]$, coprime to $\ell$ whose prime factors all lie in
$[\sqrt{z},z]$ with $\sum_{p|q} 1/p^{1-\xi/\log z}  \ge e^{\xi}/3{\xi}$ for $1\le \xi
\le \frac 23 \log z$, for which
$$
\sum_{j=1}^J \max\Sb (a,q)=1 \endSb \left|  \frac
{\vtheta((j+1)\ell/2;q,a)-\vtheta(j\ell/2;q,a) - \ell/2\vphi(q)} {\ell/2\vphi(q)}
\right|\ll  J \exp\Big(-\beta \frac{\sqrt{\log \ell}} {\sqrt{\log\log \ell}}\Big),
$$
where $\beta$ is a positive absolute constant.
\endproclaim

\demo{Proof} If $j>\epsilon (\log z)/10$ let $P_j=\emptyset$. If $1\le j\le \epsilon
(\log z)/10$ let $P_j$ be the set of primes in $I_j$, omitting the prime divisors
$v_1,\dots v_\ell$ of Siegel moduli and any prime divisors of $\ell$. Now, for $1\leq
j\leq J$ replace $x$ by $j\ell/2$ and $h$ by $\ell/2$ in Proposition 5.3, and sum,
which yields
$$
\sum_{j=1}^J \frac{1}{|{\Cal Q}|} \sum_{q\in {\Cal Q}} \max\Sb (a,q)=1 \endSb \left|
\frac {\vtheta((j+1)\ell/2;q,a)-\vtheta(j\ell/2;q,a) - \ell/2\vphi(q)}
{\ell/2\vphi(q)} \right|\ll  J \exp\Big(-\alpha \frac{\sqrt{\log \ell}} {\sqrt{\log
z}}\Big),
$$
Thus we may choose a $q$ as described in the result, proceeding as in the proof of
Corollary 5.4.
\enddemo

\proclaim{Corollary 5.6} Given a large integer $\ell$, let $Q = \ell^b$ for
sufficiently small $b>0$, and let $J$ be a positive integer $\leq \exp(\sqrt{\log
\ell})$. Select $z=10\log \ell$. Then there exists an integer $q\in [Q^{1-c/\log\log
Q}, Q]$, coprime to $\ell$ whose prime factors all lie in $[z/2,z]$ with $\sum_{p|q}
1/p^{1-\xi/\log z} \gg b e^{\xi}/\log z$ for $1\le \xi \le \frac 23 \log z$, for which
$$
\sum_{j=1}^J \max\Sb (a,q)=1 \endSb \left|  \frac
{\vtheta((j+1)\ell/2;q,a)-\vtheta(j\ell/2;q,a) - \ell/2\vphi(q)} {\ell/2\vphi(q)}
\right|\ll  J \exp\Big(-\beta \frac{\sqrt{\log \ell}} {\sqrt{\log\log \ell}}\Big),
$$
where $\beta$ is a positive absolute constant.
\endproclaim

\demo{Proof} Let $P_j=\emptyset$ for $j>1$, and let $P_1$ be a set of $\pi_1 = [2b
\log \ell/\log\log \ell]$ primes in $I_1=(z/2,z]$, omitting the prime divisors
$v_1,\dots v_\ell$ of Siegel moduli and any prime divisors of $\ell$; this is possible
since $I_1$ contains $\sim 5\log \ell/\log\log \ell$ primes, and we are forced to omit
at most $\sim \log \ell/\log\log \ell$. From here we proceed as in the proof of
Corollary 5.5.
\enddemo

\demo{Proof of Theorem 1.5} Take $Q=x^b$ in Corollary  5.4 to obtain $q$ which satisfies
the hypothesis of Corollary 3.5. Let $u:=\log y/\log z$ and select $v_{\pm}$ as in
Corollary 3.5.
 We consider the Maier matrices ${\Cal M}_{\pm}$ with
$$
({\Cal M}_{\pm})_{r,s} = \cases
\log (rq+s) & \text{\rm if} \ rq+s\ \text{\rm is prime,}\\
0 & \text{\rm otherwise,} \\
\endcases
$$
where $x/q<r<2x/q$ and $v_{\pm} \le s\le v_{\pm} +y$. Let $M_{\pm}$ denote the sum of
the entries of ${\Cal M}_{\pm}$. Using (5.3) to sum the entries in column $s$ we see
that
$$
M_{\pm}= \frac{x}{\vphi(q)} \Big( 1+O \Big(\exp\Big( -\beta \frac{\sqrt{\log
x}}{\sqrt{\log\log x}}\Big) \Big)\Big) \sum\Sb v_{\pm} \le s\le v_{\pm} + y \\
(s,q)=1\endSb 1.
$$
Let $r_+$ denote the row in  ${\Cal M}_+$ whose sum is largest, and let $x_+:= qr_+ +
v_+ \in (x,2x)$. Since there are $x/q+O(1)$ rows, we have
$$
\vtheta(x_+ + y)- \vtheta(x_+ ) \ge  \frac{q}{x} (1+O(x^{-1/2})) M_+,
$$
and then Theorem 1.5 follows from the bounds in Corollary 3.5. 
The analogous argument works for ${\Cal M}_-$.
\enddemo

\demo{Proof of Theorem 5.1} Take $Q=\ell^b$ in Corollary 5.5 to obtain $q$ which
satisfies the hypothesis of Corollary 3.3. Let $u:=\log y/\log z$ and select $U$ in
Corollary 3.3 so that $U(1-A/\log U)=u$; then put $S_{\pm}=[z^{U_{\pm}}]$, so that
$S_{\pm}\in [z^u, z^{u+Cu/\log u}]$.
 We consider the Maier matrices ${\Cal M}_{\pm}$ with
$$
({\Cal M}_{\pm})_{r,s} = \cases
\log (rq+s\ell) & \text{\rm if} \ rq+s\ell \ \text{\rm is prime,}\\
0 & \text{\rm otherwise,} \\
\endcases
$$
where $R<r<2R$ with $R:=[\ell/(2q)]$ and $1 \le s\le S_{\pm}$. Let $M_{\pm}$ denote
the sum of the entries of ${\Cal M}_{\pm}$. Using Corollary 5.5 to sum the entries in
column $s$ we see that
$$
M_{\pm} = \frac{\ell}{2\vphi(q)} \Big( 1+ O\Big( \exp\Big(-\beta \frac{\sqrt{ \log
x}}{\sqrt{\log\log x}} \Big) \Big) \Big) \sum\Sb s\le S_{\pm} \\ (s,q) =1\endSb 1.
$$
Now, the sum of the entries in row $r$ equals $0$ if $(r,\ell)>1$; and equals
$\vtheta(\ell S_{\pm} + qr; \ell, qr) - \vtheta(qr;\ell,qr) = \vtheta(\ell S_{\pm} +
qr; \ell, qr)$ if $(r,\ell)=1$, since $qr < \ell$ and $qr$ is not prime.  The number
of integers $r\in [R,2R]$ with $(r,\ell)=1$ is $R \vphi(\ell)/\ell + O(\tau(\ell)) =
\vphi(\ell)/(2q) + O(\ell^{\epsilon})$. Therefore, denoting by $r_+$ the row for which
$\vtheta(x_+; \ell, a_+)$, with $x_+=\ell S_+ + qr_+$ and $a_+=qr_+$, is maximized, we
obtain
$$
\vtheta(x_+; \ell, a_+) \geq \frac{2q}{\vphi(\ell)} (1+O(\ell^{-1/2})) M_+,
$$
and then Theorem 5.1 follows from the bounds in Corollary 3.3. The analogous argument
works for ${\Cal M}_-$.
\enddemo

\demo{Proof of Theorem 5.2} We proceed exactly as in the proof of Theorem 5.1 but
replacing the use of Corollary 3.3 by Corollary 3.4, and the use of Corollary 5.5 by
Corollary 5.6.
\enddemo

\head 6.  Further examples \endhead

\subhead 6a.  Reduced residues \endsubhead

\noindent Let $q$ be square-free.  Writing
$$
\sum\Sb n\le x\\ (n,q)=1\\ n\equiv a\pmod \ell\endSb
1 = \sum_{d|q} \mu(d) \sum\Sb n\le x\\ d|n \\ n\equiv a\pmod \ell \endSb 1,
$$
we see that this is
$$
=\cases
0 &\text{if } (a,q,\ell)>1\\
\frac{x}{\ell}\frac{\phi(q/(q,\ell))}{q/(q,\ell)} +O(\tau(q)) &\text{if } (a,q,\ell)=1.
\\
\endcases
\tag{6.1}
$$

\proclaim{Corollary 6.1} Let $q$ be a large square-free number,
which satisfies (1.3), and define $\alpha:=(\log\log q)^{-1}
\sum_{p|q} (\log p)/p$ with $\eta =\min(1/100,\alpha/3)$. Then for
$\eta (\log q)^{\eta/2} \ge u\ge 5/\eta^2$ there exist intervals
$I_\pm \subset [q/4,3q/4]$ of length at least $(\log q)^u$ such
that
$$
\sum\Sb n\in I_+\\ (n,q)=1 \endSb 1 \ge  \frac{\phi(q)}{q} |I_+|
\Big( 1+\exp\Big(-\frac u\eta(1+25\eta) \log
(2u/\eta^3)\Big)\Big),
$$
and
$$
\sum\Sb n\in I_-\\ (n,q)=1 \endSb 1\le  \frac{\phi(q)}{q} |I_-|
\Big( 1-\exp\Big(-\frac u\eta(1+25\eta) \log
(2u/\eta^3)\Big)\Big).
$$
\endproclaim

\demo{Deduction of Corollary 1.2} This follows immediately upon
noting that $z\leq (\log q)^{\eta}$ and replacing $u/\eta$ by $u$.
\enddemo

 \demo{Proof of Corollary 6.1} Take $a(n)=1$ if $n\le q$ with
$(n,q)=1$ and $a(n)=0$ otherwise.
Recall that $\eta=\min(\alpha/3,1/100)$.  Since $q$
has at most $\log q/\log \log q$ prime factors larger than $\log q$ we see that
$$
\sum\Sb p|q \\ p>\log q \endSb \frac{\log p}{p} \le 1.
$$
Therefore from our assumption that $\sum_{p|q} \log p/p =
\alpha\log \log q$ we may conclude that there exists $(\log
q)^{\eta} \le z\le (\log q)/3$ such that
$$
\sum\Sb p|q\\ z^{1-\eta}\le p\le z\endSb \frac{1}{p} \ge \eta^2. \tag{6.2}
$$

Take $\ell$ to be the product of the primes in $[z^{1-\eta},z]$
which divide $q$, so that $\ell$ is a divisor of $q$ and $\ell \le
e^{z(1+o(1))} \le q^{\frac 13}$.  Given $(\log q)^{\eta/2} \ge
u\ge 5/\eta^2$ we obtain by (6.2) and Corollary 3.2 (we check
readily that $\eta\ge z^{-1/10}$ using $\eta \ge 20\log \log \log
q/\log \log q$) that there exist points $u_\pm \in
[u,u(1+22\eta)]$ such that, with $y_\pm = [z^{u_{\pm}}]$,
$$
\sum\Sb n\le y_+ \\ (n,\ell)=1\endSb 1 \ge \Big( 1+ \exp\Big(-u(1+25\eta)
\log \Big(\frac {2u}{\eta^2}\Big)\Big) \Big) \frac{\phi(\ell)}{\ell} y_+,
\tag{6.3a}
$$
and
$$
\sum\Sb n\le y_- \\ (n,\ell)=1\endSb 1 \le \Big( 1- \exp\Big(-u(1+25\eta)
\log \Big(\frac {2u}{\eta^2}\Big)\Big) \Big) \frac{\phi(\ell)}{\ell} y_-.
\tag{6.3b}
$$

Consider now the ``Maier matrices'' ${\Cal M}_\pm$ whose $(r,s)$-th entry
is $(R+r)\ell + s$ with $1\le r\le R$ and $1\le s\le y_\pm$, and $R=[q/(4\ell)]$.
As usual we sum $a(n)$ as $n$ ranges over the elements of this matrix.  Using (6.1), note
that the $s$-th column contributes $0$ unless $(s,\ell)=1$ in which case it
contributes $R\phi(q/\ell)/(q/\ell) +O(\tau(q))$.
Thus the contribution of the matrix is
$$
\Big(R\frac{\phi(q/\ell)}{q/\ell} +O(q^{\epsilon})\Big) \sum\Sb s\le y_\pm\\
(s,\ell)=1\endSb 1.
$$
Corollary 6.1 follows immediately from (6.3 a,b).

\enddemo
\demo{Proof of Corollary 1.1} Let $q=\prod_{p\in \Cal P} p$; note
that ${\Cal A}$, the set of integers up to $q$ without any prime
factors from the set $\Cal P$, is a subset of $[1,q]$ of density
$\phi(\Cal P)/\Cal P=\phi(q)/q$, which is strictly $<1$ by
hypothesis. Let $\ell= \prod_{\sqrt{z}\le p\le z/3,\ p\in \Cal P}
p$ so that $\ell \le q^{1-\delta+o(1)}$ for some fixed $\delta>0$,
and apply the argument in our proof of Corollary 6.1 above. In
place of Corollary 3.2 we appeal to Corollary 3.3 (with $\ell$ in
place of $q$, since the hypothesis on $\Cal P$ implies that
$H(\xi)\gg \text{e}^\xi/\xi$), and see that for suitably large $u
\le \sqrt{z}$ there exist integers $y_\pm \ge z^u$ such that
$$
\sum\Sb n\le y_+ \\ (n,\ell)=1\endSb 1 \ge \Big(1+\exp(-u(\log
u+\log\log u+O(1))\Big)  \frac{\phi(\ell)}{\ell}y_+,
$$
and
$$
\sum\Sb n\le y_- \\ (n,\ell)=1\endSb 1 \le \Big(1-\exp(-u(\log
u+\log\log u+O(1))\Big)  \frac{\phi(\ell)}{\ell}y_-.
$$
We conclude that for large $u\le \sqrt{z}$ there are intervals
$I_\pm \subset [q/4,3q/4]$ of length $\ge z^u$ such that
$$
\sum\Sb n\in I_+ \\ (n,q)=1 \endSb 1 \ge \frac{\phi(q)}{q} |I_+|
\Big(1+\exp(-u(\log u+\log\log u+O(1))\Big),
$$
and
$$
\sum\Sb n\in I_-\\ (n,q)=1\endSb 1 \le \frac{\phi(q)}{q} |I_-|
\Big(1-\exp(-u(\log u+\log\log u+O(1))\Big).
$$
Finally, by Proposition 3.8 we find that we can take $y_\pm \le
z^{u+2}$, provided that $u\leq (1-\epsilon) \log\log
z/\log\log\log z$.
\enddemo

From the ``fundamental lemma'' of sieve theory (see [7]), it
follows that these estimates are essentially optimal.

\noindent {\bf Example 6.} Take $q=\prod_{p\le z} p$ and ${\Cal A}$
to be the integers below $q$ that are coprime to $q$.
We show how to tweak ${\Cal A}$ to obtain a
set ${\Cal B} \subset [1,q]$ such
that the symmetric difference $|({\Cal A}\backslash{\Cal B}) \cup
({\Cal B}\backslash {\Cal A})|$
is small, but such that ${\Cal B}$ {\sl is} well distributed
in short intervals.

Let $k=[q/(\log q)^4]$ and divide $[1,q]$ into $k$ intervals
$[mh+1,(m+1)h]$ for $1\le m\le k$, and $h= q/k =(\log q)^4 +O(1)$.
For $1\le m\le k$ consider whether
$$
\Big|\sum\Sb mh +1\le n\le (m+1)h\\ (n,q)=1\endSb 1 -
\frac{\phi(q)}{q} h \Big| \le (\log q)^3, \tag{6.4}
$$
holds or does not hold.  If (6.4) holds then take ${\Cal B}\cap [mh+1,(m+1)h]
= {\Cal A}\cap [mh+1,(m+1)h]$.  Otherwise pick an arbitrary set of
$[\phi(q)h/q]$ numbers in $[mh+1,(m+1)h]$ and take that to be
${\Cal B}\cap [mh+1,(m+1)h]$.

By construction we see that for any interval $[x,x+y]\subset [1,q]$
$$
\sum\Sb n\in {\Cal B}\\ x\le n\le x+y\endSb 1 = \frac{\phi(q)}{q}
y + O\Big( \frac{\phi(q)}q h+ \frac{y}{\log q}\Big). \tag{6.5}
$$
Thus ${\Cal B}$ is well distributed in short intervals of
length $y \ge (\log q)^5$, say.

Further note that ${\Cal A}$ and ${\Cal B}$ are quite
close to each other.  Indeed, from the Theorem in Montgomery and
Vaughan [13] we know that for integers $r\ge 1$
$$
\sum_{m\le k} \Big(\sum\Sb mh+1\le n\le (m+1)h\\ (n,q)=1\endSb 1 -
\frac{\phi(q)}{q} h\Big)^{2r} \sim \frac{(2r)!}{2^r r!}\ q
\Big(h\frac{\phi(q)}{q}\Big)^{r}. \tag{6.6}
$$
It follows that the number of values $m$ for which (6.4) does
not hold is $\ll_r q/(\log q)^{2r}$.  Therefore for any $r\ge 1$
$$
|({\Cal A}\backslash{\Cal B})\cup ({\Cal B}\backslash{\Cal A})|
\ll_r q/(\log q)^r. \tag{6.7}
$$
By (6.1) and (6.7) we therefore see that ${\Cal B}$
is also well distributed in arithmetic progressions
$a\pmod \ell$ provided $\ell \ll_r (\log q)^r$.

Now take $\ell =\prod_{\sqrt{w}\le p\le w} p =e^{(1+o(1))w}$
with $w\le z$ so that $\ell |q$.  For $u\le \sqrt{w}$ but large, our
usual Maier matrix argument then gives that one
of the following statements holds:

\noindent (i) There exists $y\in [q/4,q]$ and $a\pmod{\ell}$
 such that, with $\delta((a,\ell)=1)$ being $1$ or $0$ depending
on whether $(a,\ell)=1$ or not,
$$
\Big|{\Cal B}(y;\ell,a)- \delta((a,\ell)=1)\frac{\phi(q/\ell)}{q} y\Big|
\ge \exp(-u(\log u +\log \log u +O(1))) \frac{\phi(q/\ell)}{q} y.
$$
\noindent (ii)  There exists an interval $[x,x+y]\subset [1,q]$ with $y\ge w^u$
such that
$$
\Big|\sum\Sb x\le n\le x+y \\ n\in{\Cal B}\endSb 1 -\frac{\phi(q)}{q} y\Big|
\ge \exp(-u(\log u +\log \log u +O(1))) y.
$$
But from (6.5) we see that case (ii) cannot hold if
$w^u\ge (\log q)^5$ and if

\noindent $e^{-u(\log u+\log \log u+O(1))} \gg 1/\log q$.
We conclude therefore that the distribution of ${\Cal B}$ in
arithmetic progressions is compromised, and
that case (i) holds in this situation.  In particular
the expected asymptotic formula for ${\Cal B}(y;\ell,a)$
is false for some $\ell \ll_{\epsilon} \exp((\log q)^{\epsilon})$.

\bigskip
Our argument also places restrictions on the uniformity 
with which Montgomery and Vaughan's estimate (6.6) can hold.
Given $h$ and $q$ with $h\phi(q)/q$ large, define $\eta$ by
$q/\phi(q)=(h\phi(q)/q)^\eta$. We now show that if (6.6) holds
then $r\ll  (\log (h\phi(q)/q))^{2+4\eta+o(1)}$.

 Fix $\epsilon>0$. Choose $L$ so that $h\phi(q)/q=L^{2+\epsilon}$;
then let $u=(1-\epsilon) \log L / \log\log L$ and $w=(\log
L)^{3+4\eta+4\epsilon}$.

We follow the same argument as in example 6, but replacing
``$(\log q)^3$'' in (6.4) by ``$(h\phi(q)/q)/L$'', from which it
follows that we replace ``$y/\log q$'' in (6.5) by ``$y/L$'', and
that the upper bound in (6.7) is $\ll qh (2r/eL^\epsilon)^r$. Case
(ii) cannot hold in our range by (6.5). By the combinatorial sieve
we know that $|{\Cal A}(y;\ell,a)- \delta((a,\ell)=1)
\phi(q/\ell)y/q| \ll 2^{\pi(z)}$ so that, since case (i) holds,
$|({\Cal A}\backslash{\Cal B})\cup ({\Cal B}\backslash{\Cal A})|
\gg \phi(q)/\ell L$ in our range. By (6.7) we deduce that $
e^{w+o(w)} = \ell L (q/\phi(q)) h \gg ( eL^\epsilon/2r)^r$ so that
$r\ll w/\log L = (\log L)^{2+4\eta+4\epsilon}$, which implies the
result.

\subhead 6b.  `Wirsing Sequences' \endsubhead

\noindent Let ${\Cal P}$ be a set of primes of logarithmic
density $\alpha$ for a fixed number $\alpha \in (0,1)$; that is
$$
\sum\Sb p\le x\\ p\in {\Cal P}\endSb \frac{\log p}{p} = (\alpha+o(1))\log x,
$$
as $x\to \infty$.  Let ${\Cal A}$ be the set of integers
not divisible by any prime in ${\Cal P}$ and let $a(n)=1$
if $n\in {\Cal A}$ and $a(n)=0$ otherwise.  Wirsing proved that (see page
417 of [18])
$$
{\Cal A}(x) \sim \frac{e^{\gamma \alpha}}{\Gamma(1-\alpha)} x \prod\Sb p\le
x\\ p\in {\Cal P}\endSb \Big(1-\frac 1p\Big). \tag{6.8}
$$
Let $h$ be the multiplicative function defined by
$h(p)=0$ if $p\in {\Cal P}$ and $h(p)=1$ if $p\notin {\Cal P}$ and
take $f_q(a)=h((a,q))$ and $\gamma_q= \prod_{p|q, p\in{\Cal P}} (1-1/p)$.
Naturally we may expect that ${\Cal A}(x;q,a)
\sim \frac{f_q(a)}{q\gamma_q} {\Cal A}(x)$ and our work places restrictions
on this asymptotic.

Let $u \ge \max (e^{2/\alpha}, e^{100})$ be fixed.  Then for large
$x$ we see that the hypotheses of Theorem 2.4 are met
with $z=(\log x)/3$ and $\eta=1/\log u$.  Combining Theorem 2.4
with (6.8) which shows that ${\Cal A}(x)/x$ varies slowly (and
therefore equals $(1+o(1)){\Cal A}(y)/y$ for any $y\in (x/4,x)$)
we attain the following conclusion.

\proclaim{Corollary 6.2}  For fixed
$u \ge \max (e^{2/\alpha}, e^{100})$ and large $x$ there
exists $y\in (x/4,x)$ and an arithmetic progression $a\pmod \ell$
with $\ell \le x(3/\log x)^u$ such that
$$
\Big|{\Cal A}(y;\ell,a) -\frac{f_\ell(a)}{\ell \gamma_\ell}
{\Cal A}(y)\Big| \gg \exp(-u(\log u +O(\log\log u))) \frac{{\Cal A}(y)}
{\phi(\ell)}.
$$
\endproclaim

Similarly using Theorem 2.5 with $\eta =1/\log u$ and
$z=(\frac 13\log x)^{\frac 1M}$ we obtain the
following ``uncertainty principle'' showing
that either the distribution of ${\Cal A}$ in arithmetic progressions
with small moduli, or the distribution in short intervals must
be compromised.

\proclaim{Corollary 6.3}  Let $u \ge \max (e^{2/\alpha}, e^{100})$
be fixed and write $u=MN$ with both $M$ and $N$ at least $1$.  Then
for each large $x$ at least one of the following two statements is
true.

(i) There exists $y \in (x/4,x)$ and an arithmetic progression
$a\pmod q$ with

\noindent $q\le \exp ((\log x)^{\frac 1M})$ such that
$$
\Big|{\Cal A}(y;q,a) - \frac{f_q(a)}{q\gamma_q} {\Cal A}(y)
\Big| \gg \exp(-u (\log u +O(\log \log u))) \frac{{\Cal A}(y)}{\phi(q)}.
$$

(ii) There exists $y>(\frac 13\log x)^N$ and an interval
$(v,v+y) \subset (x/4,x)$ such that
$$
\Big|{\Cal A}(v+y)-{\Cal A}(v)-y\frac{{\Cal A}(v)}{v}\Big|
\gg  \exp(-u (\log u +O(\log \log u))) y \frac{{\Cal A}(v)}{v}.
$$
\endproclaim

\subhead 6c.  Sums of two squares and generalizations \endsubhead

\noindent{\bf Example 3, revisited.}
We return to Balog and Wooley's Example 3, the numbers
that are sums of two squares.  It is known that (see Lemma 2.1 of
[1])
$$
{\Cal A}(x;q,a) = \frac{f_q(a)}{q\gamma_q} {\Cal A}(x) \Big( 1 +
O\Big(\Big(\frac{\log 2q}{\log x}\Big)^{\frac 15}\Big)\Big). \tag{6.9}
$$
Take $q$ to be the product of primes between
$\sqrt{\log x}$ and $\log x/\log \log x$.  Using the Maier matrix
method, (6.9) and our Corollary 3.3 we obtain that for fixed
$u$ and large $x$ there exist $y_\pm \ge (\log x)^u$ and
intervals $[v_{\pm},v_{\pm}+y_{\pm}] \subset [x/4,x]$ such
that
$$
\sum_{v_+ \le n\le v_++y_+ } a(n)
\ge (1+\exp(-u(\log u+\log \log u+O(1)))) y_+ \frac{{\Cal A}(x)}{x},
$$
and
$$
\sum_{v_- \le n\le v_- + y_-  } a(n)
\le (1-\exp(-u(\log u+\log \log u+O(1)))) y_- \frac{{\Cal A}(x)}{x}.
$$
These are of essentially the same strength as the results in [1].

Further we also obtain that (for fixed $u$)
there exists $y\in (x/4,x)$ and an arithmetic progression $a\pmod \ell$
with $\ell \le x/(\log x)^u$ such that
$$
\Big| {\Cal A}(y;\ell,a)- \frac{f_\ell(a)}{\ell\gamma_\ell} {\Cal A}(y)\Big|
\ge \exp(-u(\log u+\log \log u+O(1))) \frac{{\Cal A}(x)}{\phi(\ell)}.
$$

\noindent{\bf Example 7.}\ Let $K$ be a number field with $[K:{\Bbb Q}]>1$ 
and let $R$ be its ring of integers.  Let $C_1,\dots ,C_h$ be the 
ideal classes of $R$, and define
$\Cal A^{(i)}$ to be the set of integers which are the norms of  
integral ideals belonging to $C_i$.  From the work of R.W.K. Odoni [14] 
we know that 
$$
{\Cal A}^{(i)}(x) \sim c_i \frac{x}{(\log x)^{1-E(K)}}
$$
where $c_i>0$ is a constant and $E(K)$ denotes the 
density of the set of rational primes admitting in $K$ 
at least one prime ideal divisor of residual degree $1$.  It is 
well known that $E(K)\ge 1/[K:{\Bbb Q}]$ and also we know that 
$E(K) \le 1- 1/[K:{\Bbb Q}]$ (see the charming article of J-P. Serre [16]).

We now describe  what the natural associated multiplicative functions  
$h$ and $f_q$ should be.  Define $\delta(n)=1$ is 
$n$ is the norm of some integral ideal in $K$ and $\delta(n)=0$ otherwise.
Clearly $n\in {\Cal A}^{(i)}$ for some $i$ if and only if $\delta(n)=1$. 
Naturally we would expect that 
$$
\sum\Sb n\le x\\ p^k |n \endSb \delta(n) 
\approx \frac{\sum_{j\ge k} \delta(p^j)/p^j }
{\sum_{j=0}^{\infty} \delta(p^j)/p^j} \sum_{n\le x} \delta(n),
$$
and so the natural definition of $h$ is to take 
$$
\frac{h(p^k)}{p^k}= \frac{\sum_{j\ge k} \delta(p^j)/p^j }
{\sum_{j=0}^{\infty} \delta(p^j)/p^j}.
$$
Note that $h(p)=1$ if $\delta(p)=1$ (which happens if $p$ has 
a prime ideal divisor in $K$ of residual degree $1$, and so 
occurs for a set of primes with density $E(K)$) and that 
$h(p)\le 1/p+O(1/p^2)$ 
if $\delta(p)=0$ (and this happens for a set of primes 
of density $1-E(K)>0$).  With the corresponding definition 
of $f_q(a)$ we may expect that for $(q,{\Cal S})=1$ 
(for a finite set of bad primes ${\Cal S}$ including all 
prime factors of the discriminant of $K$)  
$$
{\Cal A}^{(i)}(x;q,a) \sim \frac{f_q(a)}{q\gamma_q} {\Cal A}^{(i)}(x). 
$$
By appealing to standard facts on the zeros of zeta and $L$-functions 
over number fields one can prove such an asymptotic for 
small values of $q$ (for example, if $q$ is fixed).  Our work 
shows that that this asymptotic fails if $q$ is of size $x/(\log x)^u$ 
for any fixed $u$.  We expect that one can understand 
the asymptotics of ${\Cal A}^{(i)}(x;q,a)$ for appropriate small 
$q$ in order also to conclude that the distribution of ${\Cal A}^{(i)}$ 
in short intervals (of length $(\log x)^u$) is compromised.  
We also expect that similar results 
hold with $R$ replaced with any order in $K$.


\noindent{\bf Example 8.} Let $k$ be a fixed integer, and choose
$r$ reduced residue classes $a_1$, $\ldots$, $a_r\pmod k$ where
$1\le r< \phi(k)$.  Take ${\Cal A}$ to be the set of integers not
divisible by any prime $\equiv a_i \pmod k$ and take ${\Cal S}$ to
be the set of primes dividing $k$.  Here $h$ is completely
multiplicative with $h(p)=0$ if $p\equiv a_j\pmod k$ for some $j$
and $h(p)=1$ otherwise, and $f_q(a)$ is defined appropriately.
This is a special case of a Wirsing sequence, and so (6.8) and
Corollary 6.2 apply (we see easily that $\ell$ in Corollary 6.2
may be chosen coprime to $k$).  Note also that Example 3
essentially corresponds to the case $k=4$ and $a_1=3$.  This 
also covers Example 7 in the case when $K$ is an abelian extension.

We may apply standard techniques of analytic
number theory to study ${\Cal A}(x)$ and ${\Cal A}(x;q,a)$.  Consider
the generating function $A(s)=\sum_{n=1}^{\infty} a(n)n^{-s}$
which converges absolutely in Re$(s)>1$ and satisfies the Euler
product $\prod_{p\not\equiv a_i \pmod k} (1-p^{-s})^{-1}$.  Further
using the orthogonality relations of characters $\psi \pmod k$ we
see that
$$
A(s) = \prod_{\psi \pmod k}
L(s,\psi)^{\frac{1}{\phi(k)}\sum_{b\not\equiv a_i \pmod k}\overline{\psi(b)} } B(s),
$$
where $B$ is absolutely convergent in Re$(s)>1/2$.  Further for a character
$\chi \pmod q$ with $(q,k)=1$ we get that
$$
A(s,\chi)=\sum_{n=1}^{\infty} \frac{a(n)\chi(n)}{n^s}
= \prod_{\psi\pmod k}
L(s,\psi\chi)^{\frac{1}{\phi(k)}\sum_{b\not\equiv a_i \pmod k}\overline{\psi(b)} }
B(s,\chi),
$$
with $B(s,\chi)$ absolutely convergent in Re$(s)>1/2$.  For large $x$, if
$q\le \exp(\sqrt{\log x})$ with $(q,k)=1$ is such that
for all characters $\chi \pmod {qk}$
(primitive or not) $L(s,\chi)$ has no zeros in $\sigma \ge 1-c/\log (qk(1+|t|)$
then we may conclude by standard arguments
that
$$
{\Cal A}(x;q,a) =\frac{f_q(a)}{q\gamma_q} {\Cal A}(x) +O(x\exp(-C\sqrt{\log x})),
\tag{6.10}
$$
for some constant $1>C>0$.  Since $k$ is fixed we may suppose
that no divisor of it is a Siegel modulus.  Let $\nu_1$, $\ldots$, $\nu_t$
denote the Siegel moduli below $\exp(\sqrt{\log x})$ (see \S 5 for
details, and note that $t\ll \log \log x$), and select a prime factor
$v_i$ for each $\nu_i$.  Choose $q$ to be the product of primes between
$\sqrt{w}$ and $w$ with $w=(C/10)\sqrt{\log x}$, taking care to
omit the primes $v_1$, $\ldots$, $v_t$ which fall in this range.  Then
(6.10)  applies to this modulus $q$ (which is of size $\exp((C/10+o(1))\sqrt{\log x})$
and applying the Maier matrix method (and our Corollary 3.3) we deduce
that for large $u\le \sqrt{w}$ there exists an interval $[v,v+y]$ in $[x/4,3x/4]$
with $y\ge (\log x)^u$ such that
$$
\Big|{\Cal A}(v+y)-{\Cal A}(v) - y \frac{{\Cal A}(x)}{x} \Big|
\gg \exp(-2u(\log u+\log \log u+O(1))) y\frac{{\Cal A}(x)}{x}. \tag{6.11}
$$
Arguing more carefully, using a zero density estimate as in \S 5,
it may be possible to improve the right side of (6.11) to
$\exp(-u(\log u+\log \log u+O(1))) y\frac{{\Cal A}(x)}{x}$.

\subhead 6d. The multiplicative function $z^{\Omega(n)}$ for $z\in (0,1)$
\endsubhead

\noindent Take $a(n)=z^{\Omega(n)}$ where $\Omega(n)$ denotes
the number of prime factors of $n$ counted with multiplicity and
$z$ is a fixed number between $0$ and $1$.  We take ${\Cal S}=\emptyset$
and $h(n)=z^{\Omega(n)}$ and $f_q(a)= z^{\Omega((a,q))}$.  For
large $x$ we know from a result of A. Selberg
(see Tenenbaum [17]) that
$$
{\Cal A}(x) \sim x\frac{(\log x)^{z-1}}{\Gamma(z)}.
$$
From Theorem 2.4 and the
above,
we deduce that for fixed $u\ge \max(e^{2/(1-z)},e^{100})$
and large $x$ there exists $y\in (x/4,x)$ and an arithmetic
progression $a\pmod \ell$ with $\ell\le x(3/\log x)^u$ such
that
$$
\Big|{\Cal A}(y;\ell,a)-\frac{f_\ell(a)}{\ell\gamma_\ell} {\Cal A}(y)
\Big|
\gg \exp(-u(\log u+O(\log \log u))) \frac{{\Cal A}(y)}{\phi(\ell)}.
$$

Suppose $q\le \exp(\sqrt{\log x})$ is such that for
every character $\chi\pmod q$ (primitive or not) $L(s,\chi)$
has no zeros in $\sigma \ge 1-c/\log (q(|t|+2))$ for
some constant $c>0$.  Then following Selberg's method
we may see that for some $1>C>0$
$$
{\Cal A}(x;q,a) = \frac{f_q(a)}{q\gamma_q} {\Cal A}(x) + O(x
\exp(-C\sqrt{\log x})). \tag{6.12}
$$
Let $\nu_1$, $\ldots$, $\nu_r$ be the Siegel moduli
below $e^{\sqrt{\log x}}$ (see \S 5 for details; and note that
$r\ll \log \log x$), and select
one prime factor $v_i$ for each $\nu_i$.
Choose $q$ to be the product of primes between $\sqrt{w}$ and $w$
for $w= (C/10)\sqrt{\log x}$, taking care to
omit the primes $v_1$, $\ldots$, $v_r$ should they happen to lie
in this interval.   Then (6.12) applies to this modulus $q$,
and using the Maier matrix method and appealing to Corollary 3.3 we
find that for large $u\le \sqrt{w}$ there exists an
interval $[v,v+y]\subset [x/4,3x/4]$ with $y\ge (\log x)^u$
such that
$$
\Big|{\Cal A}(v+y)-{\Cal A}(v)-y \frac{{\Cal A}(x)}{x} \Big|
\gg \exp(-2u(\log u+\log \log u+O(1))) y\frac{{\Cal A}(x)}{x}. \tag{6.13}
$$
Taking greater care, using a zero density argument as in \S5, it
may be possible to improve the right side of (6.13) to $\gg
\exp(-u(\log u+\log \log u+O(1)))y{\Cal A}(x)/x$.

\head 7.  An uncertainty principle for integral equations \endhead

\noindent  E. Wirsing [18] observed that questions on mean-values of multiplicative
functions can be reformulated in terms of solutions to a certain integral equation.
We formalized this connection precisely in our paper [6] and we now recapitulate the
salient details.  Let $\chi: (0,\infty) \to {\Bbb C}$
be a measurable function with $\chi(t) =1$ for $0 \le t\le 1$
and $|\chi(t)|\le 1$ for all $t\ge 1$. Let
$\sigma(u)=1$ for $0\le u\le 1$ and for $u>1$ we define $\sigma$ to be the
solution to
$$
u\sigma(u) = \int_0^u \chi(t) \sigma(u-t) dt. \tag{7.1}
$$
In [6] we showed that there is a unique solution $\sigma(u)$ to (6.1) and
that $\sigma(u)$ is continuous and $|\sigma(u)| \le 1$ for all $u$.
In fact $\sigma(u)$ is given by
$$
\sigma(u) = 1+ \sum_{j=1}^{\infty} \frac{(-1)^j}{j!} I_j(u;\chi), \tag{7.2a}
$$
where
$$
I_j(u;\chi) =  \int\Sb t_1,\ldots, t_j\ge 1\\ t_1+\ldots +t_j \le u\endSb
\frac{1-\chi(t_1)}{t_1} \cdots \frac{1-\chi(t_j)}{t_j} dt_1\cdots dt_j . \tag{7.2b}
$$

The connection between multiplicative functions and the integral equation (7.1) is
given by the following result which is Proposition 1 in [6].

\proclaim{Proposition 7.1} Let $f$ be a multiplicative function with $|f(n)|\le 1$ for
all $n$, and $f(n)=1$ for $n\le y$. Let $\vartheta(x) =\sum_{p\le x} \log p$ and
define
$$
\chi(u) = \chi_f(u) = \frac{1}{\vartheta(y^u) } \sum_{p\le y^u} f(p)\log p.
$$
Then $\chi(t)$ is a measurable function with $|\chi(t)|\le 1$ for all $t$ and
$\chi(t)=1$ for $t\le 1$. Let $\sigma(u)$ be the corresponding unique solution to {\rm
(7.1)}. Then
$$
\frac{1}{y^u} \sum_{n\le y^u} f(n) =\sigma(u) +O\biggl(\frac{u}{\log
y}+\frac{1}{y^u}\biggr).
$$
\endproclaim

The converse to Proposition 7.1 also holds (see Proposition 1 (converse) of [6])
so that the study of these integral equations is entirely analogous
to the study of mean-values of multiplicative functions.  Translated
into this context our oscillation results of \S 3 take the
shape of an ``uncertainty principle'' which we will now describe.

\! We define the Laplace transform of a function $f : [0,\infty)\to {\Bbb C}$
by
$$
\Lap(f,s)=\int_0^{\infty} f(t) e^{-st} dt.
$$
If
$f$ grows at most sub-exponentially then the Laplace transform is
well-defined for complex numbers $s$ in the half-plane
Re$(s)>0$.  From equation (7.1) we obtain that for Re$(s)>0$
$$
\Lap(u\sigma(u),s) = \Lap(\chi,s) \Lap(\sigma,s). \tag{7.3}
$$
Moreover from (7.2b) we see that when Re$(s) >0$
$$
s\Lap(\sigma,s) = \exp\Big(-\Lap\Big(\frac{1-\chi(v)}{v},s\Big)\Big).
\tag{7.4}
$$
Finally, observe that if $\int_1^{\infty} |1-\chi(t)|/t\ dt<\infty$
then from (7.2b) it follows that $\lim_{u\to \infty}
\sigma(u)$ exists and equals
$$
\sigma_{\infty}:=e^{-\eta} \ \ \text{where} \ \ \eta:=\int_{1}^{\infty}
\frac{1-\chi(t)}t\ dt \ = \Lap\Big(\frac{1-\chi(v)}{v},0\Big).
$$

\proclaim{Theorem 7.2} Suppose $\sigma_\infty \neq 0$ is such that
$|\sigma(u)-\sigma_\infty| \le \exp(-(u/A)\log u)$ for some
positive $A$ and all sufficiently large $u$.  Then either $\chi(t)=1$
almost everywhere for $t\ge A$, or
$\int_0^{\infty} \frac{|1-\chi(t)|}{t} e^{Ct}dt$ diverges for some
$C\ge 0$.
\endproclaim

We view this as an ``uncertainty principle'' since (by choosing $A=1$)
we have shown that both $|\chi(t)-1|$ and $|\sigma(u)-\sigma_\infty|$
cannot be very small except in the case $\chi(t)=\sigma(u)=1$.

\demo{Proof} Since $|\sigma(u)-\sigma_\infty| \le \exp(-(u/A)\log u)$
for all large $u$ (say, for all $u\ge U$)
it follows that $\Lap (\sigma-\sigma_{\infty},s)$
is absolutely convergent for all complex $s$.
Therefore the identity
$$
s\Lap(\sigma,s) = s\Lap (\sigma-\sigma_{\infty},s) + \sigma_\infty,
$$
which {\sl a priori} holds for Re$(s)>0$, furnishes an analytic continuation
of $s{\Lap}(\sigma,s)$ for all complex $s$.
Suppose now that $\int_0^{\infty} \frac{|1-\chi(t)|}{t} e^{Ct}dt$
converges for all positive $C$.  Then $\Lap(\frac{1-\chi(v)}{v},s)$
is absolutely convergent for all $s\in {\Bbb C}$,
and so defines a holomorphic function on ${\Bbb C}$. Hence the
identity (7.4) now holds for all $s\in {\Bbb C}$.

If Re$(s)=-\xi$ then
$$
\align |s\Lap(\sigma,s)| &\le 1+|s| \int_0^{\infty} |\sigma(u)-\sigma_\infty|
e^{\xi u}du \\
&\le 1 + |s| \Big( \int_0^{U} 2e^{\xi u}du + \int_U^{\infty} \exp\Big( u
\Big(\xi-\frac{\log u}{A}\Big)\Big) du\Big)
\\
&\le 1 + |s| \Big( 2(e^{U\xi}-1)/\xi + \exp\Big(A(\xi+1)+e^{A\xi-1}/A\Big) + 1\Big),
\\
\endalign
$$
where we bounded the second integral by the sum of the two integrals
$\int_0^{e^{A(\xi+1)}}+\int_{e^{A(\xi+1)}}^{\infty}$ with the same integrand.
In the
range of the first integral one uses $u(\xi- (\log u)/A) \le e^{A\xi-1}/A$, and in the
range of the second integral one uses $u(\xi- (\log u)/A)\le -u$. Therefore, by (7.4),
if Re$(s)\geq -\xi$ and Im$(s)\ll e^\xi$ with $\xi$ large, then
$$
\text{Re } -\Lap\Big(\frac{1-\chi(v)}{v}, s\Big) \ll  e^{A\xi}.
$$

We now apply the Borel-Caratheodory lemma{\footnote {This says that for any
holomorphic function $f$ we have
$$
\max_{|z-z_0|=r} |f(z)| \le \frac{2R}{R-r} \max_{|z-z_0|=R} \text{Re} f(z)
+\frac{R+r}{R-r} |f(z_0)|
$$
where $0<r<R$.}}
 to $-\Lap(\frac{1-\chi(v)}{v},s)$
taking the circles with center $1$ and radii $r=\xi+1$ and $R=\xi+2$. Since
$$
\left|\Lap(\frac{1-\chi(v)}{v},1)\right| \leq \int_1^\infty \frac{2e^{-v}}v \le 1/2,
$$
 we deduce from the last two displayed estimates that
$$
\Big|\Lap\Big(\frac{1-\chi(v)}{v},-\xi\Big)\Big| \le \max_{|s-1|=\xi+1}
\Big|\Lap\Big(\frac{1-\chi(v)}{v},s\Big)\Big| \ll (\xi +1) e^{A\xi}.
$$
On the other hand, for any $\delta >0$ we have
$$
\Big|\Lap \Big(\frac{1-\chi(v)}{v}, -\xi\Big)\Big| \ge \int_0^\infty \frac{1-\text{Re
} \chi(v)}{v} e^{\xi v}dv \ge e^{(A+\delta)\xi} \int_{A+\delta}^{\infty}
\frac{1-\text{Re }\chi(v)}{v} dv,
$$
so that
$$
\int_{A+\delta}^{\infty} \frac{1-\text{Re }\chi(v)}{v} dv \ll \xi e^{-\delta\xi}.
$$
Taking $\delta=2\log \xi/\xi$ and letting $\xi \to \infty$, we deduce that
$\int_{A}^{\infty} \frac{1-\text{Re }\chi(v)}{v} dv =0$; that is, $\chi(v)=1$
almost everywhere for $v>A$.  This proves the Theorem.
\enddemo

\enddocument